\documentclass[12pt]{article}

\usepackage{amsmath, color}
\usepackage{amsfonts}
\usepackage{amssymb,ulem}

\newtheorem{theorem}{Theorem}

\newtheorem{definition}{Definition}
\newtheorem{lemma}{Lemma}

\newcommand{\esssup}{\mathop{\rm ess\,sup}}

\begin{document}

\def\a{\alpha}   \def\g{\gamma}  \def\b{\beta}  
\def\l{\lambda}  \def\w{\omega}   \def\W{\Omega} \def\e{\varepsilon}
\def\d{\delta}   \def\f{\varphi}  \def\D{\Delta}    \def\r{\rho}
\def\s{\sigma}   \def\G{\Gamma}   \def\L{\Lambda} \def\F{\Phi}
\def\pd{\partial} \def\ex{\exists\,}  \def\all{\forall\,}
\def\dy{\dot y}   \def\dx{\dot x}   \def\t{\tau}
\def\bx{\bar x}   \def\bu{\bar u}   \def\bw{\bar w}  \def\bz{\bar z}
\def\p{\psi} \def\bv{\bar v}   \def\k{\varkappa}
\def\calO{{\cal O}}   \def\calP{{\cal P}}  \def\calR{{\cal R}}
\def\calN{{\cal N}}   \def\calL{{\cal L}}   \def\calQ{{\cal Q}}
\def\calC{{\cal C}} \def\calM{{\cal M}} \def\calE{{\cal E}}
\def\calA{{\cal A}}
\def\calB{{\cal B}}
\def\calG{{\cal G}}  \def\calD{{\cal D}}
\def\calX{{\cal X}}  \def\calY{{\cal Y}}  \def\calZ{{\cal Z}}
\def\calF{{\cal F}} \def\calJ{{\cal J}}   \def\calV{{\cal V}}
\def\hx{\hat x}  \def\hw{\hat w}   \def\hu{\hat u}  \def\hv{\hat v}
\def\hz{\hat z} \def\1{\bf 1}  
\def\calW{{\cal W}}

\def\empty{\O}
\def\lee{\;\le\;}   \def\gee{\;\ge\;}
\def\le{\leqslant}  \def\ge {\geqslant}
\def\dis{\displaystyle}

\def\ov{\overline}  \def\ld{\ldots}
\newcommand{\dd}{\mbox{\rm\,d}}
\def\tl{\tilde}   \def\wt{\widetilde}  \def\wh{\widehat}
\def\lan{\langle}  \def\ran{\rangle}

\def\hsp{\hspace*{8mm}} \def\vs{\vskip 8mm}
\def\bs{\bigskip}     \def\ms{\medskip}   \def\ssk{\smallskip}
\def\q{\quad}  \def\qq{\qquad}
\def\np{\newpage}   \def\noi{\noindent}
\def\pol{\frac12\,}    \def\const{\mbox{const}}
\def\cl{{\rm {cl\,}}} \def\clm{{\rm {clm\,}}}  \def\supp{{\rm {supp\,}}}
\def\sec{{\rm {Sec\,}}}
\def\dist{{\rm {dist\,}}}

\def\conv{\mbox{conv\,}}
\def\vad{\vadjust{\kern-4pt}}
\def\R{\mathbb{R}}   \def\ctd{\hfill $\Box$}  \def\inter{{\rm int\,}}

\def\dis{\displaystyle}
\def\lra{\longrightarrow}  \def\iff{\Longleftrightarrow}
\def\Lra{\Longrightarrow}   \def\Lla{\Longleftarrow}
\def\iff{\Longleftrightarrow}

\def\begar{\begin{array}} \def\endar{\end{array}}
\def\be{\begin{equation}\label}  \def\ee{\end{equation}}
\def\beq{\begin{equation}\label}  \def\eeq{\end{equation}}
\def\bth{\begin{theorem}\label}   \def\eth{\end{theorem}}
\def\ble{\begin{lemma}\label}   \def\ele{\end{lemma}}
\def\bedef{\begin{definition}}    \def\endef{\end{definition}}
\def\Proof{{\it Proof.}\,}
\def\Proof{{\bf Proof.}\,}
\def\del{\sout}

\newcommand{\vmin}{\mathop{\rm vraimin}}
\newcommand{\essmin}{\mathop{\rm ess\,min}}
\newcommand{\vmax}{\mathop{\rm vraimax}}
\newcommand{\rank}{\mathop{\rm rank}}
\newcommand{\cone}{\mathop{\rm cone\,}}
\newcommand{\mes}{\mathop{\rm mes\,}}
\newcommand{\To}{\rightrightarrows}
\newcommand{\wst}{\,\stackrel{\ast}{\rightharpoonup}\,}
\newcommand{\wky}{\,\stackrel{w}{\rightharpoonup}\,}
\newcommand{\intt}{\int_{t_0}^{t_1}}

\newcommand{\blue}[1]{\textcolor{blue}{#1}}
\newcommand{\red}[1]{\textcolor{red}{#1}}

\def\bls{\baselineskip}   \def\nbls{\normalbaselineskip}

\title{\bf \large Local Minimum Principle for an Optimal Control Problem \\
with a Nonregular  Mixed Constraint}

\author{A.V. Dmitruk\thanks{Russian Academy of Sciences, Central Economics
and Mathematics Institute, Moscow, Russia;\, and Lomonosov Moscow State
University, Moscow, Russia (optcon@mail.ru).},\q
N.P. Osmolovskii\thanks{Systems Research Institute, Polish Academy
of Sciences, Warszawa, Poland (nikolai.osmolovskii@ibspan.waw.pl.). }}

\date{}
\maketitle

\begin{abstract} 
We consider the simplest optimal control problem with one nonregular mixed
constraint $G(x,u)\le0,$ i.e. when the gradient $G_u(x, u)$
can vanish on the surface $G = 0.$ Using the Dubovitskii--Milyutin theorem on
the approximate separation of convex cones, we prove a first order necessary
condition for a weak minimum in the form of the so-called ``local minimum
principle'', which is formulated in terms of functions of bounded variation,
integrable functions, and Lebesgue-Stieltjes measures, and does not use
functionals from  $(L^\infty)^*$. Two illustrative examples are given.
The work is based on the book by Milyutin \cite{AAM01}. \ssk

Keywords:\,
normed space, convex cone, dual cone, approximate separation theorem,
mixed constraint, phase point, Pontryagin function, Lebesgue-Stieltjes
measure, singular measure, costate equation.
\end{abstract}



\section{Introduction}
 Consider the optimal control problem on a fixed interval of time $[t_0,t_1]$:
\begin{eqnarray}
& {\cal J}(x,u):=\; J(x(t_0),x(t_1)) \to \min, \label{1}&\\[4pt]
\label{3} &\dot x\; =\; f(x,u),&\\[4pt]
\label{4} &G(x,u)\,\le\, 0,&
\end{eqnarray}
where the functions $J:\R^{2n}\to \R,$  $f:\R^{n+m}\to\R^n,$  and
$G:\R^{n+m}\to \R$ are continuously differentiable. This problem will
be called {\it Problem P}. \ssk

Condition (\ref{4}) is called {\it mixed state-control constraint}\, or
simply\, {\it mixed constraint}. The presence of this constraint determines
the main difficulties in obtaining necessary optimality condition for this
problem. These difficulties largely disappear if one assumes that the gradient
$G_u(x,u)$ does not vanish at the points $(x,u)\in\R^{n+m}$ where $G(x,u)=0.$
In this case we say that the mixed constraint (\ref{4}) is {\it regular}.
Traditionally, the regularity assumption {(properly modified for more
general problems)} is present in the works on necessary optimality conditions
for problems  with mixed constraints {(see e.g. \cite{MN74}--\cite{BPV}).
One of the few exceptions is the recent work \cite{Ros}, which will
be discussed later. Note that the regularity assumption for mixed constraint
does not allow one to consider the {\it pure state constraint} $g(x)\le0$
as a special case of the mixed constraint. In this paper, we do not impose
any assumptions on the mixed constraint (\ref{4}), except for the smoothness
condition for the function $G.$ \ssk

A pair $(x,u)\in \R^{n+m}$ is called {\it phase point\,} (of the mixed constraint)
if $G(x,u)=0$ and also $G_u(x,u)=0.$ As mentioned above, it is the presence of
such points, which creates the main difficulties in studying the problem and,
in addition, gives rise to the main changes even in the formulation of the
necessary optimality conditions compared to the regular case.
\ssk

We consider Problem P\, for $x\in AC([t_0,t_1],\R^n)$ and
$u\in L^{\infty}([t_0,t_1],\R^m),\;$ using the notation  \vad 
$$
\, w=(x,u)\in W\,=\; AC([t_0,t_1],\R^n)\times L^{\infty}([t_0,t_1],\R^m)
$$
and $\xi=(x_0,x_1)=(x(t_0),x(t_1)).\,$
The norm of a pair $w=(x,u)$ is 
$$
\|w\|=\; \|x\|_{AC}+\|u\|_\infty\; =\;
|x(t_0)|\,+\int_{t_0}^{t_1}|\dot x(t)|\dd t\,+\, \esssup_{t\in[t_0,t_1]}|u(t)|.
$$
Obviously, the local minimum in this norm is equivalent to the standard
{\it weak minimum}\footnote{By definition, the latter is the minimum
in the norm $||x||_C + ||u||_\infty$.}.
The goal of this paper is to obtain first-order necessary conditions for a
weak minimum in problem (\ref{1})--(\ref{4}) in the form of the so-called
{\it local minimum principle} (LMP)\footnote{Dubovitskii and Milyutin used the
term {\it local maximum principle} \cite{DM71}. Both these terms are not completely
adequate;\, nevertheless, following the authors', we use the above term.}.
\ssk

As is known, an efficient method for obtaining LMP in constrained problems
was proposed by Dubovitskii and Milyutin in \cite{DM65}. The idea was simple
(and therefore became very popular):\, at the minimum point, one should consider
the convex cones of first order approximation for the cost and constraints,
that should not intersect. Then the separation theorem is applied and the
resulting {Euler--Lagrange (stationarity) equation} is analyzed. This leads
to a LMP {with multipliers from the spaces dual to the image spaces of the
constraints}. \ssk

However, a difficulty arises in this method:\, since the image space of the
mixed constraints is $L^\infty\,,$ the separating functionals should belong
to the conjugate space $(L^\infty)^*,$ which has an essentially complex structure.
In problems  with {\it regular} mixed constraints, one can prove that the corresponding
multipliers are represented by functions from $L^1$ (see \cite{AD93, DOS14}).
Unfortunately, this is not possible for problems with {\it nonregular}
mixed constraints.

To overcome this difficulty, Dubovitskii and Milyutin \cite{DM71}\, proposed
the idea of not exact but {\it approximate} separation of the cones.
For the case of two cones, it looks as follows.\,
Let $Y$ and $X$ be normed spaces with $Y^* = X,$ let $H_0,\, H_1 \subset Y$
and $\W_0,\, \W_1\subset X$ be nonempty convex cones, $\Omega_1$ open, such
that $H_0^* =\ov \Omega_0$ and  $H_1^* = \ov \Omega_1\,,$ where the bar denotes
the closure in the strong topology of $X$ and the star denotes the dual
(conjugate) cone.

Let $x_1^0\in \Omega_1$ be a given point. Then the following is true:\, if
$\Omega_0\cap \Omega_1=\O,$ then for any $\e>0$ there exist $h_0\in H_0$ and
$h_1\in H_1$ such that $\langle x_1^0,h_1\rangle=1$ and $\|\,h_0+ h_1\|<\e.$
The converse is also true.

A similar result for a finite number of cones allowed Dubovitskii and Milyutin
to obtain in \cite{DM71} the LMP in a problem with a finite number of mixed
constraints, given as {inclusions to closed} sets in $\R^{n+m}.$
However, the book \cite{DM71} {is} published in Russian in a small number
of copies and {is} very difficult to read.
\ssk

Many years later, Milyutin presented the same result in the book \cite{AAM01},
where he considered a general problem with nonregular mixed constraints.
This time these constraints are given by smooth functions, in the form of
a finite number of inequalities like (1.3)\, and equalities $g(x,u)=0,$ assuming
that the latter satisfy the full rank condition:\, $\rank\, g_u(x,u) =\;$dim\,$g$
on the surface $g(x,u)=0,$ but without any assumpt\-ions on the joint independency
of the derivatives $G_u(x,u),\; g_u(x,u).$ The problem admits also a finite
number of endpoint constraints of the form $F(x(t_0), x(t_1))\le0$ and
$K(x(t_0), x(t_1))=0.$ Moreover, the smoothness assumption for the inequality
constraints, both the endpoint and mixed ones, were essentially weakened to
just the convexity of their directional derivatives at the reference point, while
the equality constraints were always assumed to be smooth. The problem can also
admit a pure control constraint of the inclusion type $u_2(t) \in U(t)$
on a part of control components, where the full control vector is split into two
parts: $u = (u_1, u_2),$ and $U(t)$ is a measurable set-valued mapping.
In this case, the full rank condition should be considered w.r.t. the first
group only: $\rank\, g_{u_1}(x,u_1, u_2) =\;$dim\,$g,$ as well as the
(non)regularity of all the collection of mixed constraints.
For this general problem, Milyutin obtained a necessary condition for a weak
minimum (the local maximum principle), and further developed it to
a necessary condition for a strong minimum (the global maximum principle).
A brief account of these results can be found in \cite{AD09}. \ssk 

Compared to the book \cite{DM71}, the presentation of LMP in \cite{AAM01} is
much clearer, with shorter proofs, but still difficult even for Russian-speaking
readers. Moreover, these results have never been published in English.
All this has led to the fact that the Dubovitskii--Milyutin's general theory
of the maximum principle for nonregular mixed constraints, which in our opinion
is an outstanding achievement in optimal control, still remains unknown even
to specialists.
\ssk

Because of the difficulties in the study, mentioned above, the nonregular
mixed constraints until recently remained outside the scope of specialists'
interests in the West. However, now this interest has arisen, as evidenced
by the paper \cite{Ros}. Without analyzing this publication, we will only say
that the authors did not achieve the goal that could be set: to get rid of the
functionals from $(L^\infty)^*$ in the final result, which are still present
in \cite{Ros}, though in  integral form. At the same time, the
Dubovitskii--Milyutin's LMP is devoid of this drawback. \ssk 

All this prompted us to write this article. To be as clear as possible
{in presentation of the specificity caused by the mixed constraints,
we chose the simplest possible problem for the first study:\, it includes
the Mayer cost functional, the control system, and just one mixed constraint.
(A more general problem will be considered in our future paper.)
In many ways, we follow the ideas of the book \cite{AAM01}, and yet our
presentation differs markedly from that book. We hope that this publication
will draw attention of specialists to the ideas and results contained
in \cite{AAM01}. \ssk

The paper is organized as follows. In Section~\ref{pre}, we give definitions
of the closure of a measurable set and a measurable function with respect to
the measure, proposed by Dubovitskii and Milyutin, and recall some facts
about equiintegrable sequences of functions in $L^1,$ which are used in the
proof of LMP. Section~\ref{appr} is devoted to the approximate separation
theorem for a finite number of convex cones, which plays a key role in the
proof of LMP. We formulate LMP in Section~\ref{lmp} and give two illustrative
examples in Section~\ref{ex}. The proof of LMP is given in Section~\ref{pr}.
\ms

\section{Preliminaries}\label{pre}

\subsection{The closure with respect to a measure}\,
We start with an important concept introduced by Dubovitskii and Milyutin
in \cite{DM71}. Let $M\subset\R$ be a (Lebesgue) measurable set.
The set \vad
$$
\clm M\,=\; \{t\in\R:\;\, \mes (\w\cap M)>0\q \mbox{for any open set}\;\,
\w\ni t\}
$$
is called {\it closure of $M$ with respect to (the Lebesgue) measure}.

Obviously, $\clm M$ is a closed set. Moreover, $\clm M \subset \ov M,$
but not the reverse. However, $\mes M \le \mes (\clm M)$ (since almost all
$t\in M$ are points of its density), but not the reverse.

In fact, $\clm M$ is the {\it topological support}\, of the measure $\dd\mu$
that has density $\dd\mu/\dd t = \chi_M(t),$
where $\chi_M$ is the characteristic function of the set $M$. \ssk

\if{Note also that if $M_n\,,\; n=1,2,\ld,$ is a countable family
of measurable sets, then  \vad
{$$
\bigcup_n\, \clm M_n\; \subset\; \clm \left(\bigcup_n M_n \right)\;
\subset \; \ov{\,\bigcup_n\, \clm M_n}\,,
$$
and both the inclusions are not reversible.
(The first inclusion is trivial, and the second one holds by the countable
additivity of the Lebesgue measure.) }
\ssk}\fi

Now, let $\hu: [t_0,t_1] \to\R^m$ be a measurable function.
Consider its graph $\G = \{(t,\hu(t)):\; t\in [t_0,t_1]\}$ and the projector
$ {\pi}: \R^{1+m} \to \R, \q (t,v) \mapsto t.\;$ The set
$$
\clm (\hu)\,=\; \{(t,v) \in \R^{1+m}:\;\,
\mes \pi (O\cap \G)>0\q \mbox{for any open set}\;\, O\ni (t,v) \}
$$
is called {\it the closure of the function $\hu$ with respect to (the Lebesgue)
measure}, or in short,  {\it  the closure in measure of} $\hu.$ (Note that
this definition can be applied in fact to any measurable set $\G \subset \R^{1+m},$
{and hence, to any measurable set-valued function}.)
\ssk

The following simple properties of $\clm (\hu)$ should be noted.
By $B_r(u)$ we denote the closed ball in $\R^m$ of radius $r$ centered at $u,$
and by $O_\e(t,u)$ the open set $\{ (t',u'):\; |t'-t|<\e, \;\, |u'-u|<\e\}.$
\ms

\ble{clm-pro}
If $\hu\in L^\infty([t_0,t_1],\R^m),$ then a) $\clm(\hu)$ is a compact set
in $[t_0,t_1]\times\R^m,$ which does not depend on the choice of a particular
representative of the function $\hu,$ and\, b) the projector ${\pi}$ is
surjective on $\clm(\hu),$ i.e. $\pi\, \clm(\hat u) =[t_0,t_1].$
\ele \ms

\Proof  Obviously, the set $\clm(\hu)$ is closed and bounded, which proves
the first assertion. To prove the second one, suppose the contrary, i.e.
that $\ex t_* \notin \pi\, \clm(\hat u).$ Set $r = \|\hu\|_\infty\,.$
Then, for any $u\in B_r(0)$ there is an $\e>0$ such that
$\mes \pi(O_\e(t_*,u)\cap\, \G)=0.$ Since $B_r(0)$ is compact, there exists
a finite number of $\e_i>0$ and $u_i\in B_r(0),$ $i=1,\ld, k$ such that the
union $\calC: = \bigcup_i O_{\e_i}(t_*,u_i)$ gives $\mes \pi (\calC\cap \G)=0.$
Define $\e_* = \min \e_i$ and $\w = (t_*-\e_*,\,t_*+\e_*).$
Obviously, the set $Z= \w\times B_r(0)$ is contained in $\calC,$ whence
$\mes \pi (Z \cap \G)=0.$ But the latter means that $|\hu(t)| >r$ for a.a.
$t\in \w,$ and then $\|\hu\|_\infty >r,$ a contradiction. \ctd
\ms

In fact, passing from $\hu$ to $\clm(\hu),$ we obtain a set-valued mapping
$$
\clm(\hu)(\cdot):\;\, t \mapsto\; \{v\,:\; (t,v) \in \clm(\hu)\,\}
$$
such that $\hu(t) \in \clm(\hu)(t)$ for almost all $t\in [t_0,t_1].$
Clearly, this mapping is upper semicontinuous.

Another way to define the closure in measure is
$
\clm (\hu):= \bigcap_{u\sim \hat u}\, \ov{ \mbox{Graph}\, (u)},
$
where the equivalence $u\sim \hat u\,$ means that $u(t)= \hu(t)$ for almost
all  $t\in  [t_0,t_1].$ One can easily show that  \vadjust{\kern-8pt}
$$
\clm (\hu)\,= \bigcap_{\mes E=\,t_1-t_0}\, \ov{ \mbox{Graph}\, (\hat u|_E)}\,,
$$
where this time the intersection is taken over all measurable
subsets $E\subset [t_0,t_1]$ of  full measure, and by definition
$\mbox{Graph}\, (\hat u|_E)=\{ (t,u)\in\R^{1+m}:\;\, t\in E,\;\, u=\hat u(t)\}.$
\ssk

Note also that, if $\hx(t)$ is a continuous function, then  \vad
\beq{xclmut}
\clm (\hx,\hu)(t)\, =\, (\hx(t),\clm (\hu)(t))\q\;
\mbox{for all}\q t\in [t_0,t_1].
\eeq

\subsection{Uniformly integrable families of functions}\, A family $\calF$
of functions from $L^1([t_0,t_1],\R^m)$ is called {\it uniformly integrable}
(or {\it equi-integrable}) if for any $\e>0$ there is a $\d>0$ such
that for any measurable set $E \subset [t_0,t_1]$ of $\mes E<\d$ we have
$\int_E |\,\l(t)|\dd t<\e,\q \all \l\in \calF.$
\ssk

Obviously, this is equivalent to the fact that the functions of $\calF$
possess a common modulus of integrability, i.e. a function
$\nu:\, \R_+ \to \R_+$ such that $\nu(\d)\to0$ as $\d\to 0+,$ and for any
measurable set $E \subset [t_0,t_1]$ we have
$$
\int_E |\,\l(t)|\dd t \lee \nu(\mes E),\qq \all \l\in \calF.
$$
\ssk

Since the functions of $ L^1([t_0,t_1],\R^m)$ generate absolutely continuous
vector-valued measures on $[t_0,t_1],$ a uniformly integrable family of functions
generates a uniformly absolutely continuous family of vector-valued measures.
\ssk

We will use these concepts in the case when $\calF$ is a sequence of functions
$\l^k\in L^1([t_0,t_1],\R^m),$ $k=1,2,\ld.$ By the Dunford--Pettis theorem
\cite{DS, Ed}, any uniformly integrable sequence $\l^k\in L^1$ contains
an $L^\infty$-weakly convergent subsequence $\l^{k_s}.$
The latter means that there exists a function $\l\in L^1$ such that,
for any $ u\in L^\infty$ we have \vad
$$
\int_{t_0}^{t_1}\lan \l^{k_s},u \ran\dd t\;\, \to\;
\int_{t_0}^{t_1}\lan \l,u \ran\dd t \q\; (s\to\infty).
$$
We write in this case $\; \l^{k_s}\wky \l \q\; (s\to\infty).$
\ssk

In general,  it is impossible to extract a weakly convergent sequence from
an arbitrary bounded set of functions in $L^1,$ since this set can be not
uniformly integrable. Nevertheless, the following important fact holds true
(see, e.g. \cite{Saad-Val} and references therein)\footnote{\,Dubovitskii
and Milyutin, being not aware of these works, proved this fact independently
in  \cite{DM71, AAM01}.}.
\ssk

\ble{bit}(The biting lemma.)\,
Let a sequence $\l^k \in L^1([t_0,t_1],\R)$ be bounded, i.e. $\|\l^k\|_1\le\const$
for all $k=1,2,\ld.$ Then there exists a sequence of measurable sets
$A^k\subset [t_0,t_1]$ such that $\mes A^k\to (t_1-t_0)$ and the sequence
$\l_A^k:= \l^k\chi_{A^k}$ is uniformly integrable, hence it contains a weakly
convergent subsequence.
\ele
\ssk

\subsection{Functions of bounded variation and charges}\,
Denote by $\R^{n*}$ the space of row vectors of the dimension $n,$ and by
$BV([t_0,t_1],\R^{n*})$ the space of  functions $p:[t_0,t_1]\to \R^{n*}$ of
bounded variation for which the values $p(t_0-0)$ and $p(t_1+0)$ are also
defined. By definition, {\it the jump} of $p$ at a point $t\in[t_0,t_1]$ is
$[p](t):=p(t+0)-p(t-0).$ We define the Radon measure (or charge) $\dd p,$
which corresponds to the function $p,$ by the following condition:\,
if $[a,b]\subset [t_0,t_1],$ then
$
\int_{[a,b]} \dd p\; =\; p(b+0)\,-\,p(a-0).
$
Note that  we always prefer to denote Radon measures on $[t_0,t_1]$
by $\dd p,$ \if{$\dd \mu,$}\fi rather than $p$ \if{$\mu$}\fi or $p(\dd t),$
\if{$\mu(\dd t),$\fi as is customary.
This makes it possible to distinguish measures from the functions of
bounded variation that define them, without introducing new notation.
\ssk

If $l: C([t_0,t_1],\R^{n})\to \R$ is a linear continuous functional, then
by the Riesz theorem, there exists a function $p\in BV([t_0,t_1],\R^{n*})$
such that  \vad
\beq{int}
\langle l,x\rangle =\int_{t_0}^{t_1} x(t)\dd p\qq \forall \,
x\in C([t_0,t_1],\R^{n}),
\eeq
but this function is not unique\footnote{\,Strictly speaking, we should write
$\dd p\,x,$ but it is more convenient to write $x \dd p.$}.
It is unique under the additional requirement that the function $p$ vanishes
at $t_0$ (or at $t_1$) and is one-way continuous, for example continuous from
the left. For the definiteness, we will assume that the functions $p\in BV$
are left-continuous, i.e., $p(t-0)=p(t)$ for all $t\in [t_0,t_1],$ and $p(t_0)=0.$
If $p$ belongs to $BV,$ we write $\dd p\in C^*,$  keeping in mind the
relations (\ref{int}).
\ssk

If the function $p\in BV([t_0,t_1],\R^{n*})$ is absolutely continuous, then
the measure $\dd p$ and the functional $l$ defined by (\ref{int})
are also called {\it absolutely continuous}.
\ssk

We say that $t_*\in[t_0,t_1]$ is a {\it point of continuity} of a measure
$\dd p$ if $[p](t_*)=0,$ i.e., if the measure $\dd p$ has no atom at this
point. Recall that the  measure $\dd p$ can have atoms in at most countably
many points. Therefore, the set of continuity points of $\dd p$ is dense in
$[t_0,t_1]$. A point $t^*\in[t_0,t_1]$ where the  measure $\dd p$ has an atom,
i.e., $[p](t_*)> 0,$ is often called a {\it jump point} of the measure.
\ssk

As usual, we say  that a sequence of measures $\dd p^k$ weakly* converges
to a measure $\dd p\in C^*$ (i.e. $C-$converges in $C^*$)\, if
$$
\int_{t_0}^{t_1} x(t)\dd p^k\, \to\, \int_{t_0}^{t_1} x(t)\dd p\q\; \mbox{as}
\q k\to\infty
$$
for all $x\in C([t_0,t_1],\R^{n}),$ and we write in this case
$\dd p^k \wst \dd p.$ \ssk

Let $\l^k$ be a sequence of functions in $L^1$. Consider the corresponding
sequence of absolutely continuous measures $\dd p^k:= \l^k\dd t.$
Assume that $\dd p^k$ is weakly* convergent to some measure $\dd p\in C^*,$
that is $\l^k\dd t \wst \dd p$. Denote by $\Theta\subset [t_0,t_1]$ the set
of all continuity points of the measure $\dd p$. Then for any
$\t_0,\t_1\in \Theta$ with  $\t_0<\t_1\,,$ we have  \vad
$$
\int_{[\t_0,\t_1]}\l^k(t)\dd t\;\to\, \int_{[\t_0,\t_1]}\dd p \q (k\to\infty).
$$ \ssk

\section{An approximate separation theorem}\label{appr}
Let $X$ and $Y$ be normed spaces, such that $X=Y^*.$ Let $\Omega\subset X$
be a nonempty convex cone and $\overline\Omega$  its closure.  We say that
a cone $H\subset Y$ is {\it thick on the cone $\Omega$} (or is {\it predual to}
$\W)$ if $ H^*=\, \ov\Omega.\,$ Here $H^*$ denotes the {conjugate} cone
of $H,$ consisting of all linear continuous functionals that are nonnegative
on $H$ (in other words, $ H^*$ is the polar cone of $(-H)$).
We will need the following properties of these cones.  \ssk

\ble{sec}
Let $x^0 \in \inter \W.$ Then the set $Sec\,H =\, \{ h\,:\; \langle x^0,h\rangle=1 \}$
is bounded, and its conical hull is $H \setminus\{0\}.$
\ele \ssk

\Proof\, Suppose $\ex h_k \in Sec\,H$ with $\|h_k\| = r_k \to\infty.$
Setting $\wt h_k= h_k/r_k$ we have $\|\wt h_k\|=1$ and $\lan x^0, \wt h_k\ran\to 0.$
Let $x^0+B_\d \subset\W$ for some $\d>0,$ where $B_\d$ is the closed ball in $X$
of radius $\d,$ centered at zero. Then $\langle x^0+B_\d,\wt h_k\rangle \ge 0,$
whence $\langle B_\d, \wt h_k\rangle \ge -o(1).$ But here the infimum of the
left hand side equals $ -\d,$ a contradiction.

Thus, $Sec\,H$ is bounded. To prove the second assertion, take any nonzero
$h\in H.$ Since $H^*= \ov \W,$ we have $\langle \W, h \rangle \ge 0,$
and since $x^0 \in \inter \W,$ we have $\a:= \langle x^0, h \rangle > 0,$
whence $h/\alpha\in \sec H,$ q.e.d. \ctd
\ssk

Now, let be given two convex cones $H_0\,,\, H_1\subset Y$ and two convex cones
$\W_0\,,\, \W_1 \subset X,$ such that $H^*_0 = \ov\W_0$ and $H^*_1 = \ov\W_1\,,$
where again $X=Y^*.$ The following theorem is an approximate analog of the
Hahn--Banach separation theorem for the case of two convex cones, in which the
separating functionals are taken not from the dual but from the predual space.
\ssk

\bth{septwo}\,
Let $\W_1$ be open and $x^0_1\in\W_1\,.$ Then $\W_0\cap\, \W_1 =\O$
$\iff$ $\all\e>0$ $\ex (h_0,h_1)\in H_0\times H_1$
such that $\langle x^0_1,h_1\rangle =1$ and $\|h_0+h_1\|<\e.$
\eth \ssk

\Proof $(\Lla)$ Suppose, on the contrary, that $\ex\hx\in \W_0\cap\, \W_1\,.$
Without loss of generality assume that $\|x^0_1\|= \|\hx\|=1$ and
$B_r(\hx)\subset\W_1$ for some $r>0,$ where $B_r(\hat x)$ is the closed ball
in $X$ of radius $r$ centered at $\hat x.$ Set $\e= r/2$ and take any pair
$(h_0,h_1)$ with the above properties. They imply $\|h_1\|\ge 1.$
Set $y= h_0+h_1\,.$ Then $\|y\|\le\e$ and $h_0+h_1-y=0,$ whence \vad
\beq{hh01}
\langle \hx,h_0\rangle\, + \langle\hx,\, h_1-y\rangle\,=\,0.
\eeq
The first summand here is nonnegative. Now, the inequality
$\langle\hx-B_r, h_1 \rangle\ge 0$ implies that
$\langle\hx, h_1\rangle\ge \sup \langle B_r, h_1\rangle \ge r,$
and since $|\langle\hx, y\rangle| \le \e,$ the second summand in (\ref{hh01})
can be estimated as $\lan \hx,\, h_1-y \ran \ge r -\e = r/2 >0,$ so the sum
in (\ref{hh01}) cannot be zero, a contradiction.

$(\Lra)$ We have to show that $\,\inf \|H_0+ Sec\,H_1\|=0.$ Suppose the
contrary:\, this $\inf >r>0,$ i.e. the distance from the set $Sec\,H_1$
to the cone $-H_0$ is greater than $r.$ Therefore,  
$
(-H_0)\,\cap\, (Sec\,H_1 +B_r)\,=\, \O.
$
Then, by the classical separation theorem, $\ex\hx\in X,\; \|\hx\|=1,$
such that $\lan \hx, -H_0 \ran \le0$ and $\lan\hx,\,Sec\,H_1 + B_r\ran\ge 0.$
The first relation implies $\hx\in H_0^* = \ov\W_0\,,$ and the second one
$\lan\hx,\,Sec\,H_1 \ran \ge \sup \lan\hx,B_r\ran =r,$ i.e.
$\all h_1\in Sec\,H_1$ we have $\lan \hx, h_1\ran \ge r.$ By Lemma~\ref{sec}
$\|Sec\,H_1\| \le d$ for some $d>0.$ Take any positive $\e<r/d\,.$
Then it follows that for any $h_1\in Sec\,H_1$  \vad
$$
\lan \hx +B_\e\,,h_1\ran \gee r -\, \sup\, \lan B_\e\,,h_1\ran \gee r-\e d >0,
$$
and since the conical hull of $Sec\,H_1$ is $H_1\setminus\{0\},$
we get $\hx+B_\e \subset H_1^* =\, \ov\W_1\,,$ so $\hx \in \W_1.$
Thus, $\hx \in \ov\W_0 \cap \W_1\,,$ and since $\W_1$ is open, there
exists an element $x' \in \W_0 \cap \W_1\,,$ a contradiction.
\ctd \ms

Note that in the proof of implication $\Lra,$ instead of separating the
cones $\W_0$ and $\W_1$ by an element of $X^*,$ we use the classical theorem
to separate the cone $H_0$ and an extension of the cone $H_1$ by an element
of $X.$ \ssk
\ms

{\bf The general case.}\, Now, let be given a finite number of convex
cones $\W_0\,, \W_1,$ $\ldots,$ $\W_m$ in $X,$ among which the last $m$
are open, and convex cones $H_0$,$H_1,\,\ld,\,H_m$ in $Y$ such that
$H_i^* = \ov\W_i$ for all $i=0,1,\,\ld,\,m$ (i.e. each $H_i$ is thick on $\W_i).$
As before,  $X= Y^*.$ Let be also given elements $x^0_i \in\W_i\,,$
$i=1,\,\ld,\,m,$ of the open cones.\, The following theorem is an approximate
analog of the Dubovitskii--Milyutin ``multi-separation'' theorem for convex
cones\footnote{{In \cite{DM65}, the condition for separating the cones was
called the Euler--Lagrange equation.}}
(see \cite[Theorem 2.1]{DM65}).
\ssk

\bth{sepmany}
$\W_0\cap\, \W_1\cap\,\ld\,\cap\,\W_m =\O$ 
$\iff$ $\all\e>0\;\; \ex h_i\in H_i\,,$ $i=0,1,\,\ld,\,m,$ such that \vad
\beq{sum1e}
\sum_{i=1}^m\langle\, x^0_i,h_i\rangle\,=1 \qq \mbox{and} \qq \|\,h_0+ \sum_{i=1}^m h_i\|<\e.
\eeq
\eth

The first of these conditions can be regarded as a {\it normalization condition},
while the second one is an {\it approximate Euler--Lagrange equation}.
Note that the cone $\W_0$ does not appear in the first condition,
it appears only in the second one. \ssk

\Proof $(\Lla)$ Without loss of generality assume that $\sum_1^m \|x^0_i\|= 1.$
We have to show that all $\W_i$ do not intersect. Suppose, on the contrary,
$\ex\hx\in \bigcap_{i=0}^n\, \W_i\,.$ Without loss of generality  assume that $\|\hx\|=1,$ and
let  $r>0$ be such that $B_r(\hx)\subset\W_i$ for all $i\ge 1.$

Set $\e= r/2\,$ and take any collection $(h_0, h_1,\ldots,h_m)$ {of elements
in $H_0,H_1,\ldots,H_m$, respectively,} satisfying
(\ref{sum1e}). The first of these conditions together with the relation
$\sum_1^m \|x^0_i\|= 1$ imply that $\dis\,\max_{1\le i\le m}\|h_i\|\ge 1.$
Let, for definiteness, $\|h_m\|\ge1.$
\ssk

Set $y=\sum_0^m h_i\,.$ Then  $\|y\|\le\e$ and $\sum_0^{m-1} h_i +  h_m -y  =0,$
whence \vad
\beq{hh0m}
\sum_{i=0}^{m-1}\langle\; \hx, h_i\rangle\, +\, \langle\hx,\, h_m-y\rangle\,=\,0.
\eeq
The first $m$ summands here are nonnegative.
Now, the inequality $\langle\hx-B_r, h_m \rangle\ge0$ implies that
$\langle\hx, h_m \rangle\ge \sup \langle B_r, h_m\rangle\ge r,$ and since
$|\langle\hx,y\rangle| \le \e,$ the last summand in (\ref{hh0m}) can be
estimated as $\langle\hx,\, h_m-y\rangle \ge r -\e = r/2>0,$ so
the left hand side in (\ref{hh0m}) cannot be zero, a contradiction.

$(\Lra)$ We prove by induction. Suppose the theorem holds
for all $m'<m$ open cones and consider the case of $m$ open cones.

If $\W_0\cap\, \W_m =\O,$ then by Theorem \ref{septwo} $\all \e>0$
$\ex h_0\in H_0$ and $h_m\in H_m$ such that $\langle x^0_m,h_m\rangle =1$ and $\|h_0+h_m\|<\e.$
Choosing all $h_i$ for $i=1,\,\ld,\,m-1$ to be arbitrary sufficiently small
elements of $H_i\,,$ we get $\sum_1^m\langle x^0_i,h_i\rangle \ge 1$ and
{$\|h_0+ \sum_1^m h_i\|< 2\e.$}
Obviously, this implies the statement of Theorem.
\ssk

Now, suppose that $W_0:= \W_0\cap\, \W_m \ne \O.$ Set $K_0= H_0 + H_m$
and notice that in this case
$K_0^* = \ov\W_0 \cap\, \ov\W_m = \ov{\W_0 \cap \W_m} = \ov W_0\,,$ that is
$K_0$  is thick on $W_0.$ (The second equality holds because both $\W_0$
and $\W_m$ are convex and the last one is open.)\,
Consider the cones $K_0\,,H_1\,,\ld,\, H_{m-1}\subset Y$ and the corresponding
cones $W_0\,,\W_1\,,\ld,\, \W_{m-1}\subset X,$ {where the last collection
does not intersect}. The dual cones to the first ones are equal to the closure
of the last ones, so we have the situation of Theorem \ref{sepmany}\, for $m-1$
open cones. By the premise of induction, $\all\e>0$ $\ex k_0\,, h_1\,,\ld,\, h_{m-1}$
from the cones $K_0\,, H_1\,,\ld,\, H_{m-1}\,,$ respectively, such that
$$
\sum_1^{m-1}\langle x^0_i,h_i\rangle =1 \q \mbox{and }\q 
\|\,k_0+ \sum_1^{m-1} h_i\|<\e.
$$
Setting here $k_0= h_0 + h_m$ with some $h_0\in H_0$ and $h_m\in H_m\,,$ 
we obtain
$
\sum_1^{m-1}\langle x^0_i,h_i\rangle + \langle x^0_m, h_m\rangle \ge 1 $ 
and still  $\|h_0+ \sum_1^{m-1} h_i + h_m\| < \e.
$
Multiplying the obtained collection by some $\l\le1$ we get the required. \ctd
\ms

\section{Local minimum principle} \label{lmp}
Consider the set  \vad
$$
{\calN(G)}\,:=\; \{(x,u)\in\R^{n+m}:\q G(x,u)=0, \q G_u(x,u)=0 \}.
$$
Clearly, $\calN(G)$ is closed.\, It is called  the {\it set of phase points}.
We assume that {this set} is nonempty (otherwise the mixed constraint is
regular).
\ssk

Define {the following set-valued} mapping
$(x,u)\in\R^{n+m}\,\To\, S(x,u)\subset \R^{n*}:$ \ssk

\begin{itemize}
   \item[(i)] if    $(x,u)\in \calN(G),\;$ then $S(x,u)=\,\{G_x(x,u)\},$
   \item[(ii)] if $(x,u)\notin \calN(G),\;$ then $S(x,u)=\,\O.$
\end{itemize}  \ssk

\noi
Thus, $S(x,u)$ is  a singleton $\{G_x(x,u)\}$ or an emptyset.
\ssk

For any nonempty set $M\subset\R^{n+m}$ we define
$
S(M)\,=\; \bigcup_{(x,u)\in M} S(x,u).$

Let $\hat w= (\hat x,\hat u)\in W$ be a given admissible process in Problem P
investigated for optimality. Denote for short $\hat \xi=(\hat x(t_0),\hat x(t_1)).$
Let us formulate the conditions of local minimum principle for the process $\hat w.$
\ssk

Recall that for the function $\hat u$ we introduced (in Sec. \ref{pre})
the set-valued mapping
$
\clm(\hat u)(t)\,=\; \{u\in\R^m:\; (t,u)\in \clm(\hat u)\}, 
$
and recall also that
$\big(\hat x(t),\clm(\hat u)(t)\big)\,=\; \clm(\hat w)(t)$
for all $t\in[t_0,t_1].\,$ Define a set  \vad
\beq{setD}
\calD\,:=\; \{t\in[t_0,t_1]:\q
\clm (\hw)(t)\, \cap\, \calN(G)\ne\O\}.
\eeq
Since the set $\clm (\hw)$ is compact and $\calN(G)$ is closed, $\calD$ is
a closed (possibly empty) subset in $[t_0,t_1].$
Denote by $\chi_\calD$ its characteristic function.
\ssk

For any $t\in\calD$ consider the set
$\conv S(\clm(\hat w)(t))= \conv S\big(\hx(t),\clm(\hu)(t)\big),$ where
$\conv$ stands for the convex hull. We call it the {\it  set of possible
directions of jumps}\, of the adjoint variable at the point $t.$
\ssk

For any {nonempty} set $M\subset\R^{n+m}$ we define
$
G_x(M)\,=\; \bigcup_{(x,u)\in M} G_x(x,u).$ \\
It follows from the definitions that for any $t\in\calD$ we have
$$ S(\clm(\hat w)(t))\;=\;   G_x\big(\clm(\hat w)(t)\cap \calN(G)\big) \ne\O.
$$

Now, define {\it the Pontryagin function} $\; H(x,u,p)\; =\; p\,f(x,u),\,$
where $p\in \R^{n*}$ is a costate (adjoint) row-vector. \ssk

The conditions of {\it local minimum principle} (LMP) at the point $\hat w$
are as follows:\, there exist multipliers  \vad
\be{15}\qq
\hat \alpha_0\in\R, \q \hat p\in BV([t_0,t_1],\R^{n*}), \q
\hat \l\in L^1([t_0,t_1],\R),
\q \dd\hat \eta\in C^*([t_0,t_1],\R), \q
\ee
such that  \vad
\be{16} \hat \alpha_0\ge0,\q \hat \l\ge0,\q \hat \l\, G(\hat w)=0,\q
\dd\hat\eta\ge0, \q \dd\hat\eta\cdot\chi_\calD\,=\dd\hat\eta,
\ee
\be{21}
\hat \a_0\,+\, \|\hat \l\|_1\,+\int_{[t_0,t_1]}\dd\hat \eta\;>\,0,
\ee
{and a $\dd\hat\eta$-measurable essentially bounded function
$\hat s: [t_0,t_1]\to \R^{n*}$ such that
\be{17b}
\hat s(t)\,\in\, \conv S\big(\hx(t),\clm(\hu)(t)\big) \q
\mbox{for almost all}\;\; t \;\; \mbox{in}\;\, d\hat\eta-\mbox{measure},
\ee }
moreover, the following adjoint equation in terms of measures:
\be{18}
-\dd \hat p\,=\; H_x(\hat w,\hat p)\dd t\,+\,
\hat \l\, G_x(\hat w)\dd t\, +\,\hat s\dd\hat \eta,
\ee
and the transversality conditions:  \vad 
\be{19}
-\hat p(t_0-)\,=\; \hat \a_0 J_{x_0}(\hat \xi), \qq
\hat p(t_1+)\,=\; \hat \a_0 J_{x_1}(\hat \xi)
\ee
are fulfilled, and finally, the stationarity condition with respect
to the control is satisfied:  \vad
\be{20}
H_u(\hat w,\hat p)\,+\,\hat \l\, G_u(\hat w)\,=\,0.
\ee

The last equation means that
$
H_u(\hat w(t),\hat p(t))\,+\,\hat \l(t)\,G_u(\hat w(t))\,=\,0,\q
\mbox{a.e. in} \q [t_0,t_1],
$
{where} "a.e." means "almost everywhere with respect
to the Lebesgue measure".
\ssk

Condition (\ref{17b}) means that there exists a set $\calR \subset \calD$
of full $\dd\hat\eta$-measure (i.e.,
$\int_\calR \dd\hat\eta =\int_{[t_0,t_1]}\dd\hat\eta$)\, such that
the inclusion in (\ref{17b}) holds for all $t\in \calR$. \ssk

{The values $\hat s(t)$ for $t\notin \calD$ are of no importance.}
\ssk

Note that conditions (\ref{15})--(\ref{20}) differ from that for problems
with {\it regular}\, mixed constraints only by the presence of the term
$\hat s\dd\hat \eta\,$ in the adjoint equation (\ref{18}).
If $\calD=\O$, this term vanishes.

The adjoint equation can be understood in the integral form:\,
for almost all $t$  \vad
$$
\hat p(t)\,=\; \hat p(t_0-0)\,+ \int_{t_0}^t \bigl(H_x(\hat w,\hat p)+
\hat \l\, G_x(\hat w)\bigr)\dd \t\, +\int_{t_0-0}^{t+0} \hat s(\t)\,\dd\hat \eta(\t).
$$
(The last integral is taken over the interval $[t_0,t]$ including its endpoints.)
\ms

{\it Remark.}\, It is convenient 
to introduce the so-called {\it augmented Pontryagin function}
$\q \ov H(x,u,p,\l)\; =\; p\,f(x,u)\, +\, \l\, G(x,u),\;$
whence the costate equation (\ref{18}) and the stationarity condition in
the control (\ref{20}) take the following shorter form, respectively: \vad
\begin{eqnarray}\label{18b}&&
-\dd \hat p\,=\; \ov H_x(\hat w,\hat p,\hat\l)\dd t\,+\,\hat s\dd\hat \eta,\\[4pt]
&& \label{20b}
\q\ov H_u(\hat w,\hat p,\hat\l)\,=\,0.
\end{eqnarray}\ssk

\begin{theorem}\label{thmp}
If $\hat w$ is a weak local minimum in Problem P, then it satisfies the local
minimum principle (\ref{15})--(\ref{20}).
\end{theorem}
\ms

{\bf Some particular cases.} \ssk

1. Let $t_*$ be an isolated point in $\calD$ and the function $\hu$ be
continuous at $t_*\,.$ \newline
Then $\clm(\hu)(t_*) = \hu(t_*)$ and $s(t_*) = G_x(\hw(t_*)),$
so the measure $\dd\hat\eta$ can have an atom at this point:
$\dd\hat\eta(\{t_*\})>0,$ and the costate
variable {have} the jump $[p](t_*) = - G_x(\hw(t_*))\,\dd\hat\eta(\{t_*\}).$
\ssk

If the control $\hu$ has a discontinuity of the first kind at $t_*\,,$ then
$(\clm\hu)(t_*)$ consists of two points: $\hu(t_*-0)$ and $\hu(t_*+0).$
{If both the corresponding points $(\hx(t_*),\hu(t_*-0))$ and $(\hx(t_*),\hu(t_*+0))$
belong to $\calN(G),$ then }
$s(t_*) = s_0\,G_x(\hx(t_*),\hu(t_*-0)) + s_1\,G_x(\hx(t_*),\hu(t_*+0)),$
where $s_0\ge0,\; s_1\ge0,$ $ s_0+s_1=1,$ and the costate variable has
the jump $[p](t_*) =\,- s(t_*)\,\dd\hat\eta(\{t_*\}).$
\ms

2. Consider the case when the function $G$ does not depend on $u,$
i.e. $G(x,u) = g(x).$ Then the mixed constraint (\ref{4}) reduces to
{\it a pure state constraint} $g(x)\le 0.$ In this case
$\calN(g) =\, \{(x,u):\; g(x)=0\},$ i.e., each point on the boundary of the
state constraint\footnote{\,To be precise, this is indeed the boundary
if $g'(x)\ne 0$ on {it}. } is a phase point, the set
$\calD =\, \{t:\, g(\hx(t))=0\}$  {consists of
contact points}, and by setting $\calR=\calD,$ we get $s(t) = g'(\hx(t))$
at any point $t\in \calD.$
\ssk

Consequently, the formulation of LMP in this case is as follows:\, there exist
multipliers $\hat \alpha_0\in\R,$ $\hat p\in BV([t_0,t_1],\R^{n*}),$
$\hat \l\in L^1([t_0,t_1],\R),$ and $\dd\hat \eta\in C^*([t_0,t_1],\R)$
such that \vad
$$\begar{rcl}&&
\hat \alpha_0\ge0,\q\, \hat \l\ge0,\q\, \hat \l\, g(\hat x)=0,\q\,
\dd\hat\eta\ge0, \q\, g(\hat x)\dd\hat\eta=0,
\\[4pt]
&& \hat \a_0\,+\, \|\hat \l\|_1\,+\int_{[t_0,t_1]}\dd\hat \eta\;>\,0,
\\[6pt]
&& -\dd \hat p\,=\; H_x(\hat w,\hat p)\dd t+
 g'(\hat x) \big(\hat \l\dd t+\dd\hat \eta\,\big),
\\[4pt]
&& -\hat p(t_0-)=\, \hat \a_0 J_{x_0}(\hat \xi), \q\;
\hat p(t_1+)=\, \hat \a_0 J_{x_1}(\hat \xi),
\\[4pt]
&&  H_u(\hat w,\hat p)\,=\,0,
\endar
$$
where $H(x,u,p)=pf(x,u).$
Setting  $\hat \l\dd t+\dd\hat \eta =:\dd\hat\mu,$  we get \vad
$$\begar{rcl}&&
\hat \alpha_0\ge0,\q\dd\hat\mu\in  C^*([t_0,t_1],\R),\q
\dd\hat\mu\ge0, \q g(\hat x)\dd\hat\mu=0, \q
\hat \a_0+\int_{[t_0,t_1]}\dd\hat \mu>0,
\\[4pt]
&&  -\dd \hat p\,=\; H_x(\hat w,\hat p)\dd t\,+\, g'(\hat x) \dd\hat\mu.
\endar
$$
The transversality conditions and the condition $H_u(\hat w,\hat p)=0$
do not change.  Thus we obtain the well-known conditions of LMP for the
problem with a state constraint.
\ms

\section{Examples}\label{ex}

\subsection*{Example 1:\, the measure $\dd\hat\eta$ has atoms}
Let $[t_0,t_1]$ be a fixed interval, $t_0<t_1,$ $\,x\in\R,$ $u\in\R.$
Consider the problem \vad
$$  J:=\; x(t_0)\,x(t_1) \to \min, \qq
\dot x=u,\qq G:=\; \frac{1}{2}\,u^2 -x +1\le0.
$$
{Conditions $G=0,\;\, G_u=0$ select here the only phase point }
$(x,u)=(1,0).$ Since $(x-1) \ge \pol u^2,$ we always have $x\ge1,$
hence $\inf J\ge1.$ Then, the process $\hx(t)\equiv 1,$ $\hu(t)\equiv 0$
is a solution to the problem. Therefore, $\calD=[t_0,t_1].$ \ssk

Further, we have (removing the hats over the multipliers):
$
H=\, p\,u,$  $ \ov H =\, p\,u\, +\l\, G,$ $ \ov H_u=\, p+ \l\, u,$
$ \ov H_x=\,-\l,$ $ s=G_x= -1.$
The condition $\ov H_u =0$ gives $p+\l\hat u=0, $ whence $p(t)=0$
for all $t\in(t_0,t_1),$ and therefore $p(t_0+) =\, p(t_1-)=\,0.$
The transversality conditions
give $p(t_0-)=\,-\a_0\, x(t_1)=\,-\a_0,$ $ p(t_1+) =\, \a_0\,x(t_0)=\, \a_0,$
so the jumps of $p$ at the endpoints are:
$[p](t_0):=\, p(t_0+)-p(t_0-)=\,\a_0$ and $[p](t_1):=\,p(t_1+)-p(t_1-) =\, \a_0.$
\ssk

The adjoint equation $-\dd p= \ov H_x\dd t +s\dd\eta$ reduces to
$\dd p= \l \dd t + \dd\eta,$ $ \l\ge0,$ $\dd\eta\ge 0.$
Since  $p(t)=0$ for $t\in(t_0,t_1),$ we have $\l=0$ and $\dd p= \dd \eta.$
\ssk

If $\a_0=0,$ then $\dd p=\dd \eta=0,$ which contradicts the nontriviality
condition (\ref{21}). 
Therefore, we can set $\a_0=1.$ Then the measure $\dd p\,$ is the sum of
$\d$-functions at $t_0$ and $t_1,$ respectively, and the same is true
for $\dd\eta.$
\ms

\subsection*{Example 2:\, the measure $\dd\hat\eta$ is absolutely continuous}

Fix any $T>0$ and consider the problem on the interval $[-T, T]$: \vad
\beq{examp2} \begar{l}
\dis \dy=\,x, \q\; \dx=\,u,\qq G(y,x,u)=\, \pol u^2 -x \le 0, \\[8pt]
\dis  J\,=\; y(T) -y(-T)\, -\, \frac m 2 \Bigl( x(-T) + x(T) \Bigr) \;\to\;\min,
\endar
\eeq
where $m\in (0,T)$ is a given number.

Here $y(T) -y(-T) = \int_{-T}^T x\,dt,$ so the variable $y$ is in fact
inessential. The set $\calN(G) = \{ (x,u):\; G=0,\;\, G_u=0\}$ consists of the
only point $(x,u)= (0,0),$ and since $(G_y, G_x) = (0,-1),$ the direction
of possible jumps of the costate vector $p = (p_y, p_x)$
is $s = (s_y, s_x) = (0, -1),$ where the subscripts indicate coordinates,
not partial derivatives. \ssk

Set $b= T-m,$ and consider the following process: \\[4pt]
\hsp $\hx(t) = \hu(t) = 0\;$ on the interval $[-b,b],$ \\[4pt]
\hsp $\hx(t) = \pol (t-b)^2$ and $\hu(t) = t-b\;$ on $[b,T],$ \\[4pt]
\hsp $\hx(t) = \pol (t+b)^2$ and $\hu(t) = t+b\;$ on $[-T,-b].$ \ms

Obviously, this process is admissible. Let us show that it is globally
optimal in the problem. To do this, choose any value $h\in {[0,T]},$
fix the endpoints $x(-T) = x(T)= \pol h^2,$ and find the minimum of
$\int_{-T}^T\, x\,dt\,$ under the given mixed constraint $G\le 0,\;$ i.e.
$|u|\le \sqrt{2\,x}.$ Clearly, this minimum is attained at the lowest
possible curve, i.e. the one satisfying $\,\dx = \sqrt{2\,x}\,$ on $[0,T],$
$x(T)=h,\,$ and symmetrically on $[-T,0].$ Therefore,  \vad
\beq{x-sym}
\begar{l}
x(t)=\, 0\q \mbox{on the interval}\;\; [-T+h,\,T-h], \\[4pt]
x(t)=\, \pol(t-T+h)^2\,\q \mbox{on}\;\; [T-h,\,T], \\[4pt]
x(t)=\, \pol(-t+T-h)^2\,\q \mbox{on}\;\; [-T,\,-T+h].
\endar
\eeq

\noi  Then $J(h) =\, {\frac 13 h^3 - \frac m 2\, h^2},$  and we have to find
the minimum of this function over $h\in {[0,T]}.$ The equation
$J'(h) =\, h^2 - m h =0$ has the only {positive} solution $h=m,$ and since
$J'(h)<0$ for $h<m,$ and $J'(h)>0$ for $h>m,$ we conclude that $J(h)$ has a
global minimum over $h\in [0,T]$ at $h=m.$ Clearly, no $h>T$ can give a smaller
cost value, so $h=m$ provides the global minimum of $J(h),$ and the
corresponding curve (\ref{x-sym}), coinciding with $\hx(t),$ provides
the global minimum in problem (\ref{examp2}). \ssk

Let us check the LMP for this curve.
According to Theorem \ref{thmp}, there exist $\a_0\ge0,$ $\l\in L^1,$
$\l(t)\ge0$ a.e. on $[-T,T],$ a measure $\dd\eta \in C^*$ supported on
$\calD= [-b,b],$ the function $s = (s_y, s_x) = (0,-1)$ a.e. in $[-b,b]$
w.r.t.$\,\dd\eta,$ and the function $p = (p_y, p_x) \in BV,$ such that
this collection is nontrivial:  \vad
\be{6b}
\a_0\,+\, \|\l\|_1\, +\int_{[t_0,t_1]}\dd\eta\, >0,
\ee
generates the augmented Pontryagin function
$\ov H =\, p_y x+ p_x u + \l\, (\pol u^2 -x),$
and satisfies the conditions (\ref{18})--(\ref{20}).
\ssk

The condition $\ov H_u=0$ gives
$ p_x\,+\,\l u\,=\,0.\; $
Since $\ov H_y=0$ and $s_y=G_y=0,$ the adjoint equation
$-\dd p_y=\, \ov H_y\dd t +s_y\dd\eta\,$ reduces to $\dd p_y=0,$
whence $p_y=\const.$ The transversality conditions for $p_y$ are:
$p_y(-T)=\,-\a_0 J_{y(-T)}=\a_0$ and $p_y(T)=\,\a_0 J_{y(T)}=\a_0.$
Consequently, $p_y=\,\a_0.$ \ssk 

Since $H_x=p_y=\,\a_0$ and $s_x=G_x= -1$ (a.e. in $[-b,b]$ w.r.t. $\dd\eta),$
the adjoint equation $-\dd p_x=\, \ov H_x\dd t +s_x\dd\eta\,$ has the form
\be{1ap}
-\dd p_x\, =\, \a_0\dd t - \l \dd t -\dd\eta.
\ee
The transversality conditions for $p_x$ are: \vad
\be{trp}
p_x(-T)=\, -\a_0 J_{x(-T)}=\, \a_0\,m/2\,,\qq
p_x(T)=\, \a_0 J_{x(T)} =\, -\a_0\,m/2\,.
\ee

If $\a_0=0,$ then $p_y \equiv 0,$ and (\ref{1ap}) with (\ref{trp}) reduce
to $p_x = \l\dd t + \dd\eta \ge 0$ and $p_x(-T)= p_x(T)=0,$ whence
$\l=0$ and $\dd\eta=0,$ which contradicts the nontriviality (\ref{6b}).
Thus, we can set $\a_0=1,$ and also $p_y =1.$ Then
$$
\dd p_x\,=\, (\l-1)\,\dd t + \dd\eta, \qq p_x(-T) = m/2\,, \q\; p_x(T) = -m/2\,.
$$
Since $u=0$ on $D=[-b,b]$ and $p_x = -\l u,$ we get
$p_x=0$ and $\l\dd t + \dd\eta =\dd t$ there. So, $\dd\eta$ is absolutely
continuous on $D$ and is not unique:\, both $\l$ and $\dot\eta$ are just
nonnegative and bounded by the relation $\l(t) + \dot\eta(t) =1.$
\ssk

Consider the interval $(b,T].\;$ We have there $\dd\eta=0,$ \vad
$$
x =\, (t-b)^2/2,\q\; u =\dx =\, t-b,\q\; p_x =\, -\l\, u =\, -\l (t-b),\q\; \dot p_x = \l-1,
$$
whence $\dot\l (t-b)=\, 1-2\l.$ Setting $\s = \l- 1/2$ and $\t= t-b,$
we get $\dot\s\, \t =\, -2\s,$ which easily gives $\s =\, c/\t^2,$ and so,
$\l =\,\pol +\, c/\t^2$ with some constant $c.$ Then $c=0$ and $\l = 1/2$
(otherwise $\l \notin L^1),$ whence $p_x =\,-\pol (t-b)<0$ and $p_x(b+0) =0,$
so the jumps $[p_x](b) =[\eta](b)=0.$
\ssk

The symmetric picture is on the interval $[-T,-b).\;$ Here $\l = 1/2,$
$p_x = \,-\pol (t+b)>0$ and $p_x(-b-0) =0,$ so the jumps
$[p_x](-b) =[\eta](-b)=0.$ 
\ms

\section{Proof of LMP} \label{pr}
In this section we prove Theorem \ref{thmp}. We will assume that
\be{a} \esssup_{t\in[t_0,t_1]}\, G(\hat w(t))\,=\,0,
\ee
otherwise the mixed constraint is redundant  for the weak minimality of
the process $\hw.\;$ For any $\d>0,$ define a set  \vad
$$
M_\d\,=\; \{t\in[t_0,t_1]:\;\; G(\hat w(t))\ge-\d\}.
$$
In view of (\ref{a}), $\;\mes M_\d>0$ for all $\d>0.$
\ms

\subsection{Application of approximate separation theorem}

Let us consider as independent variables in Problem P the pair
$(x_0, u)\in \R^n\times L^\infty,$ while the state $x(t)$ is determined
by the latter as the solution to equation (\ref{3}) with the initial
condition $x(t_0)= x_0\,,$ so that $x = x(x_0,u)$ is a nonlinear
operator of $(x_0,u),$ which maps $\R^n\times L^\infty$ {to the space $C.$ }
The Problem P has then the form
\beq{ProbP}
J(x_0,\, x(x_0,u)(t_1)) \to \min, \qq G(x(x_0,u)(t),\, u(t)) \le 0.
\eeq
Note that the weak minimality of the pair $\hw =(\hx,\hu)$ in Problem P
is equivalent to the  {local} minimality of the pair
$(\hx(t_0),\hu)$ in Problem (\ref{ProbP}).
\ssk

1. Let $\hw=(\hx,\hu)$ be a reference process. Consider the equation in
variations for the control system (\ref{3}):  \vad
\beq{eqvar}
\dot{\bx}\, =\; f_x(\hat w)\,\bx\, +\, f_u(\hat w)\,\bu,
\qq \bx(t_0)=\bar x_0\,,
\eeq
and define the corresponding linear operator
$$
A:\; (\bar{x}_0,\bu)\in \R^n\times L^\infty([t_0,t_1],\R^m)\;
\to\; \bx\in C([t_0,t_1],\R^n),
$$
where $\bx$ is the solution to (\ref{eqvar}) for the given pair
$(\bar{x}_0,\bu).$
\ssk

Recall the following well-known fact, which relates to the  nonlinear
operator $(x_0,u) \to x$ defined by the original equation (\ref{3})
with $x(t_0)=x_0\,.$ This operator maps a neighborhood of the point
$(\hat x(t_0),\hat u)\in \R^n\times L^\infty([t_0,t_1],\R^m)$ to the space
$C([t_0,t_1],\R^m)$ endowed with its standard norm $\|x\|_C=\, \max_t |x(t)|$.
\ssk

\ble{frechet}
The operator $A$ is the Frechet derivative at $(\hat x(t_0),\hat u)$
of the nonlinear operator $(x_0,u) \to x.$
Hence, for any solution $\bw =(\bx,\bu)$ to (\ref{eqvar}),
there is a correction $\wt x_\e$ parametrized by $\e>0$ with
$\wt x_\e(t_0)=0$ and $\|\wt x_\e\|_C = o(\e)$ as $\e\to 0+,$ such that
the pair $w_\e = (\hx+ \e \bx+ \wt x_\e,\, \hu +\e \bu)$ satisfies
(\ref{3}) with the initial condition $\hat x_0 +\e \bar x_0\,.$
\ele
\ms

2. Introduce a Banach space
$\calY= \R^{n*}\times L^1([t_0,t_1],\R^{m*})\times\R$ with elements
$y=(c_0,v,r),$ and its dual space
$\calX= \calY^* =\R^n\times L^\infty([t_0,t_1],\R^m)\times\R$
with elements $\k=(x_0,u,q).\;$ The pairing between these spaces
is given by
$$
\lan y,\k\ran\,=\; c_0\,x_0\, +
\int_{t_0}^{t_1} v(t)\,u(t)\,\dd t\,+\,r q.
$$

To prove Theorem \ref{thmp}, we follow the Dubovitskii--Milyutin approach.
First of all, we define, in the space $\calX,$ the following cones of
first order approximations of the cost and constraint.\,
For any $\d>0$ we set  \vad
$$ \begar{ccl}  \Omega_0 & = & \{\k\in\calX :\;\;
J_{x_0}(\hat \xi)\,x_0\, + J_{x_1}(\hat \xi)\,x(t_1)\,+ q<0,
\;\; \mbox{where}\;\;  x=A(x_0, u)\,\}, \\[6pt]
\Omega_\d & = & \dis \{\k \in \calX:\;\;
\esssup_{t\in M_\d}\,\Bigl(G_x(\hat w)\,x +\,G_u(\hat w)\,u\Bigr)\,+ q\, <0,
\;\; \mbox{where}\;\;  x=A(x_0, u)\,\},
\\[10pt]
\Omega & = & \{\k\in\calX:\;\;  q >0\}.
\endar
$$

Obviously, all these cones are convex, open, and nonempty (since the first two
contain the triple $(0,0,-1),$ and the last one contains $(0,0,1)$).\,
\ms

{\small
{\it Remark.}\, In what follows, our aim will be to separate these cones by
elements of the predual space $\calY.$ The variable $q$ and the third cone $\W$
are introduced because without $q$ the second cone $\W_\d$ can be empty, which
prevents application of the separation theorem. To avoid the analysis of this
case that can be tedious, we, following \cite{AAM01}, introduce the additional
variable $q>0.$ The price for this trick is negligible in the case of present
simplest problem~P. In a more general problem, it would be a bit more
essential, but still acceptable.}
\ms

The first step in the Dubovitskii--Milyutin approach is to show that the
cones of first order approximations do not intersect. \ssk

\begin{lemma}\label{th5}
If $\hat w$ is a weak minimum in problem~P, then for any $\d>0$
\be{22}
\Omega_0\,\cap\, \Omega_\d\, \cap\, \Omega\; =\; \O.
\ee
\end{lemma}

\Proof\, Suppose, on the contrary, there exist a $\d>0$ and a \vad
$$
\bar\k\,=\; (\bar x_0, \bu, \bar q)\,\in\, \Omega_0\cap \Omega_\d\cap\Omega.
$$

Set $\bx=A(\bar x_0,\bu),$ $\bw=(\bx,\bu),$
and take the curve $w_\e = (\hx+ \e \bx+ \wt x_\e,\, \hu +\e \bu)$ from
Lemma \ref{frechet}. Since $G(\hw)\le0$ a.e. on $[t_0,t_1],$ and
$G'(\hw)\bw =\, G'_x(\hw)\,\bx + G'_u(\hw)\,\bu < -\bar q<0$ on $M_\d\,,$
we have for sufficiently small $\e>0$:
$$
G(w_\e) =\; G(\hw) + \e G'(\hw)\bw + o(\e)\; <\, -\e\bar q + o(\e) <0 \q
\mbox{a.e. on}\; M_\d\,.
$$
For a.a. $t\notin M_\d\,,$ we have $G(\hw)\le -\d,$ whence we obviously
obtain $G(w_\e)\le -\d/2 <0$ for small $\e>0,$ so the pair $w_\e$ satisfy the
mixed constraint of problem (\ref{ProbP}).
\ssk

Now, consider the reference endpoints $\hat\xi = (\hat x_0,\, \hat x_1)$
and set $\bar x_1 = \bx(t_1),$ $\bar\xi = (\bar x_0,\, \bar x_1).$
Since $\bar\k \in \W_0$ and $\bar q>0,$ we have $J'(\hat\xi)\bar\xi =
J_{x_0}(\hat \xi)\bar x_0 + J_{x_1}(\hat \xi)\bar x_1 < -\bar q <0,$
and then, for sufficiently small $\e>0$  \vad
$$
\begar{c}
\calJ(w_\e) =\;
J(\hat x_0 +\e\bar x_0,\, \hat x_1 + \e\bar x_1 +\wt x_\e(t_1))\;= \\[6pt]
=\; J(\hat\xi) + \e J'(\hat\xi)\,\bar\xi\, +
\, J'_{x_1}(\hat\xi)\,\wt x_\e(t_1)+ o(\e)\;\;
<\; J(\hat\xi) - \e\,\bar q + o(\e)\; <\; \calJ(\hw),
\endar
$$
which contradicts the weak minimality at $\hw.$ The lemma is proved. \ctd
\ssk

3. Next, we define cones $H_0,\, H_\d,\,H$ in $\calY$ that are thick on
$\Omega_0,\, \Omega_\d,\,\Omega,$ respectively. \,
Let us start with the cone $\W_0\,.$  Consider a functional
$l:\R^n\times L^\infty \to\R$ such that  \vad
$$
l( x_0, u):= \, J_{x_1}(\hat \xi)\, x(t_1), \q
\mbox{where}\q x=A( x_0, u).
$$
As is known, one can give its explicit dependence of $(x_0, u).$
To this end, introduce the usual adjoint function $p_0\in AC$ determined
by the adjoint equation to (\ref{eqvar}):
$$
-\dot p_0\,=\,p_0f_x(\hat w) \q \mbox{with}\q p_0(t_1)= J_{x_1}(\hat \xi).
$$
Then, obviously $\frac d{dt}(p_0\,x)\,=\, p_0 f_u(\hw)\,u,$ whence integrating
we get
$$
J_{x_1}(\hat \xi)\, x(t_1)\,  =\; p_0(t_0)\,x_0\,+
\int_{t_0}^{t_1}p_0\,f_u(\hat w)\, u\dd t
\qq \forall\, ( x_0, u)\in \R^n\times L^\infty.
$$
Consequently, \vad
\beq{ome-0}\q
J_{x_0}(\hat \xi) x_0+ J_{x_1}(\hat \xi)\,x(t_1) +q\,=\,
\Bigl(J_{x_0}(\hat \xi) + p_0(t_0)\Bigr)x_0\,+
\int_{t_0}^{t_1}p_0f_u(\hat w)\, u\dd t\,+q
\eeq
for all $ \k=( x_0, u,q )\in \calX.\;$ Define a triple  \vad
$$
\hat y_0\,=\, (J_{x_0}(\hat \xi) + p_0(t_0),\;
p_0 f_u(\hat w),\;1 )\, \in\, \calY.
$$
In view of (\ref{ome-0}), $\Omega_0$ is an open half-space:
$\Omega_0=\{\k\in\calX: \;\lan \hat y_0,\k\ran<0\},$ and its closure is\,
$\ov \Omega_0=\{\k\in\calX: \;\lan \hat y_0,\k\ran\le0\}.$
Setting $\; H_0\,=\,\{-\a_0\,\hat y_0:\;\, \a_0\ge0 \},\;$
we obtain $H_0^*=\, \ov \Omega_0\,,\;$ that is $H_0$ is thick on $\Omega_0$.
\ssk

4. Consider the cone $\W_\d\,.$ First, we claim that  
$$
\ov \Omega_\d\,=\; \{\k\in\calX:\;\;
\esssup_{t\in M_\d}\,\Bigl(G_x(\hat w)\, x +\,G_u(\hat w)\,u\Bigr)\,+  q\, \le 0,
\;\; \mbox{where}\;\;  x=A(x_0, u)\, \}.
$$
Indeed, for any such $\k,$ taking a smaller $q'<q\,$ we get a point
$\k' \in \W_\d\,,\;$ q.e.d.
\ssk

Define a cone $ H_\d$ consisting of all functionals $y_\d = (c_0, v, r)\in\calY$
that for all $ \k=(x_0, u, q)\in \calX$ act as follows:
$$
\langle y_\d,\k\rangle\, =\; -\int_{t_0}^{t_1}
{\l(t)}\Bigl(G_x(\hat w)\, x\, + \,G_u(\hat w)\,u\, +q \Bigr)\,dt,\q
\mbox{where} \q x= A(x_0,u),
$$
and $\l\in L^1$ is an arbitrary nonnegative function concentrated on $M_\d,$
that is $\l\ge0$ and $\l\chi_{M_\d}=\l,$ where $\chi_{M_\d}$ is the
characteristic function of the set ${M_\d}$.
\ms

\begin{lemma}\label{th6}
    $\; H_\d^*\; =\; \ov  \Omega_\d$.
\end{lemma} \ssk

\Proof\,  If $ y_\d\in H_\d$ and $\k\in\ov \Omega_\d\,,$ then obviously
$\langle y_\d,\k\rangle\ge 0,$ whence $\k\in H_\d^*.$
Therefore, $\ov \Omega_\d\subset  H_\d^*\,.$ Let us prove the converse inclusion
$\ov \Omega_\d\supset  H_\d^*\,.$
\ssk

Indeed, take any $\bar\k\in H_\d^*,$ that is $\lan y_\d,\bar\k\rangle\ge0$
for all $ y_\d\in H_\d.$ This means that
$$
\langle y_\d,\bar\k\rangle\, =\; -\int_{t_0}^{t_1} {\l(t)}
\Bigl(G_x(\hat w)\,\bx\, + \,G_u(\hat w)\,\bu\, +\bar q \Bigr)\,dt \gee 0
$$
for all nonnegative functions $\l\in L^1$ concentrated on $M_\d\,.$
This obviously implies
$\;
G_x(\hat w)\, \bx\, + \,G_u(\hat w)\,\bu\, +\bar q \,\le 0\,
$
a.e. on $M_\d\,,$ that is $\bar\k\in\ov \Omega_\d\,.\;$
Thus, $H_\d^*\subset \ov \Omega_\d\,,\;$ q.e.d. \ctd
\ssk

5. Take any $y_\d\in H_\d$ and the corresponding function $\l\in L^1.$
Represent it in the canonical form $y_\d = (c_0, v, r).$ In fact, we only
have to find a representation of the term $\int \l(t)\, G_x(\hw)\,x\dd t.$
To this aim, define a function $p_\d \in AC$ from the equation
$$
-\dot p_\d\,=\; p_\d f_x(\hat w)\,+\,\l\, G_x(\hat w), \q p_\d(t_1)=0.
$$
Since $\dx = f_x(\hw)\,x + f_u(\hw)\, u,$ we have
$
\frac {\dd}{\dd t}(p_\d\,x)\,=\, -
\l\, G_x(\hat w)\, x\,+\, p_\d f_u(\hw)\,u,
$
whence
$$ \intt \l\, G_x(\hat w)\, x\,\dd t\; =\; p_\d(t_0)\, x_0\, +
\intt p_\d f_u(\hw)\, u\,\dd t.
$$
Then, for any $\k=(x_0,u,q)\in \calX$  we have \\[4pt]
$\dis \langle y_\d,\k\rangle =\; -p_\d(t_0)\, x_0\, -
\intt \Big(\big(p_\d f_u(\hat w) +
\l G_u(\hat w ) \big)\, u\,+ \l q \Big) \,\dd t, $
and so
$$
y_\d =\; -\left(p_\d(t_0),\;\; p_\d f_u(\hat w) +\l\, G_u(\hat w ),\;\;
\int_{t_0}^{t_1}\l \dd t\right).
$$
\ssk

6. Finally,  consider the cone $\W.$ Set $\hat y =(0,0,1)\in \calY$
and $H=\{\a\hat y: \a\ge 0 \}.$ Then  $ H^*=\ov \Omega,$ that is $H$
is thick on $ \Omega.$
\ssk

7. Set $ \k^0=(0,0,1)\in \calX$ (here $x_0=0,$ $u=0,$ $q=1$).
Obviously,
$$
-\k^0\in \, \Omega_0\,\cap\, \Omega_\d\,, \qq \k^0 \in \W.
$$

Now, we apply Theorem \ref{sepmany} to condition (\ref{22}).
According to this theorem, for any $\d>0$ and any $\e>0$ there
exist functionals
\be{23}
y_0\in  H_0,\q y_\d\in  H_\d,\q y\in H
\ee
such that
\begin{eqnarray}&&\label{24}
\lan y_0,\, -\k^0\ran\, +\,  \lan y_\d,\, -\k^0\ran\; =\,1,
\\[4pt]
&& \label{25} \|y_0\,+\,y_\d\,+\,y\|\,<\e.
\end{eqnarray}
(Here we choose the cone $\W$ to be excluded from the normalization condition
(\ref{24}).) \ssk

Analysis of these conditions will lead to the local minimum principle.

\subsection{Analysis of conditions (\ref{23})--(\ref{25})}

According to the definitions of $H_0,$ $H_\d,$ and $H,$
conditions (\ref{23}) mean that
$$
\begar{ccl}
y_0 &=& -\a_0\big(p_0(t_0)+ J_{x_0}(\hat \xi),\;
p_0f_u(\hat w),\; 1 \big),\q \a_0\ge0,\\[6pt]
y_\d &=& -\big(p_\d(t_0),\;\;p_\d f_u(\hat w) + \l\, G_u(\hat w ),\;\;
\int_{t_0}^{t_1}\l \dd t\big),\q \l\ge0,\q \l\chi_\d= \l,\\[6pt]
y &=& \a\,(0,\;0,\;1),\q \a\ge0.
\endar
$$
Condition (\ref{24}) gives \vad
\be{26}
\a_0\,+\int_{t_0}^{t_1}\l \dd t\;=\,1.
\ee
In view of this relation, we get \vad
$$
\begar{c}
-(y_0+y_\d+y)\, =\, \\[4pt]
=\; \Big(\a_0p_0(t_0)+p_\d(t_0) +\a_0J_{x_0}(\hat \xi),\;\;
(\a_0p_0+p_\d )f_u(\hat w)+ \l G_u(\hat w ),\;\; 1-\a \Big).
\endar
$$
Set $p=\a_0p_0+p_\d.$  Then
\be{27a}
-\dot p\,=\, pf_x(\hat w)+ \l G_x(\hat w),
\qq p(t_1)=\a_0J_{x_1}(\hat\xi),
\ee
and 
$-(y_0 +y_\d +y)\, =\, \Big(p(t_0)+\a_0J_{x_0}(\hat\xi),\;\,
pf_u(\hat w)+ \l G_u(\hat w ),\;\, 1-\a \Big).$ \\
Condition (\ref{25}) implies
\be{28}
|\,p(t_0)+\a_0J_{x_0}(\hat\xi) |\,+\,
\|\,pf_u(\hat w)+ \l G_u(\hat w )\|_1 \,+\, |1-\a|\;<\,\e.
\ee
Recall that such $\,\a_0\,,\, \a,\, \l,\, p\;$ exist for all $\d>0$
and $\e>0.$
\ssk

Thus, there exist two countable sequences
$
\{(\a_{0}^k,  \a^k, \l^k, p^k)\}_{k=1}^\infty $
{and} $\{\d^k\}_{k=1}^\infty,
$
where $\d^k\to 0+,$
\begin{eqnarray}\label{29}&&
\a_{0}^k\in\R,\q \a^k\in\R,\q \l^k\in L^1, \q p^k\in AC,
\\[4pt]
&&\label{30}
\a_{0}^k\ge0,  \q \a^k\ge0, \q \l^k\ge 0,\q \l^k\chi_{M_{\d^k}}=\l^k,
\\[4pt]
&& \label{Mdk}
-\d^k\le G(\hat w(t))\le 0 \q \mbox{a.e.\, on}\;\; M_{\d^k}\,
\end{eqnarray}
(the latter follows from the definition of $M_{\d}),$ such that
\begin{eqnarray}\label{31}&&
\a^k\to 1, \q\;  \a_{0}^k\,+\|\l^k\|_1\, =\,1,\\[4pt]
&& \label{27}
-\dot p^k\,=\; p^k f_x(\hat w)\,+\, \l^kG_x(\hat w),
\\[4pt]
&& \label{32}
p^k(t_0)+ \a_{0}^kJ_{x_0}(\hat\xi)\to0, \qq
p^k(t_1)=\, \a_{0}^kJ_{x_1}(\hat\xi),
\\[4pt]
&& \label{32a} \|\,p^kf_u(\hat w)\,+\, \l^kG_u(\hat w )\|_1\,\to 0.
\end{eqnarray}

Hereinafter, we {do not write} the condition $k\to\infty.$
Note also that superscript $k$ is always used to denote the number of a member
in the sequence and never used to denote the degree.
\ssk

Without loss of generality we assume that {$\;\a_0^k\, \to\, \hat \a_0\gee 0.$}
Then
\be{34a}
{\hat\a_0\, +\, \|\l^k\|_1 \, \to\, 1.}
\ee
Moreover, conditions (\ref{32}) imply
\be{tr}
p^k(t_0)\to\, -\hat\a_0J_{x_0}(\hat\xi), \qq
p^k(t_1)\to\, \hat\a_0J_{x_1}(\hat\xi).
\ee
It follows that the sequences $p^k(t_0)$ and $p^k(t_1)$ are bounded, and
in view of (\ref{34a}) the norms $\|\l^k\|_1$ are also bounded. Therefore,
by (\ref{27}) and the Gronwall's inequality, the norms $\|p^k\|_\infty$
are uniformly bounded as well.  \bs

Now, we rewrite the adjoint equation (\ref{27}) in the form of measures:
\be{35}
-\dd p^k\,=\; p^kf_x(\hat w)\dd t\, +\, \l^k G_x(\hat w)\dd t.
\ee
Define a measure  \vad
\be{dmk}
\dd\mu^k\,:=\; \l^k\, G_x(\hat w)\dd t.
\ee
Equation (\ref{35}) then takes the form
\be{36}
-\dd p^k\,=\; p^kf_x(\hat w)\dd t\, +\, \dd\mu^k.
\ee
Clearly, the sequence $\|\dd\mu^k\|$ is bounded. Without loss of generality
we assume that $\dd\mu^k$ weakly* converges to some measure $\dd\hat\mu\in C^*$
(i.e. $C-$converges in $C^*$), and denote this as  \vad
\be{35a}
\dd\mu^k \wst \dd\hat\mu.
\ee

Conditions (\ref{35a}), (\ref{36}), and (\ref{tr}) imply 
that there is a function $\hat p\in BV$ such that at every point $t\in[t_0,t_1]$
of continuity of the limiting measure $\dd\hat\mu$ (hence almost everywhere)
we have $p^k(t)\to \hat p(t),$ and moreover,
\begin{eqnarray}\label{38}&&
-\dd \hat p\; =\; \hat p\,f_x(\hat w)\dd t\,+\, \dd\hat \mu,
\\[4pt]
&& \label{39}
-\hat p(t_0-)\, =\, \hat\a_0J_{x_0}(\hat\xi), \qq
\hat p(t_1+)\,=\, \hat\a_0J_{x_1}(\hat\xi).
\end{eqnarray}
Since the sequence $\|p^k\|_\infty$ is bounded, we also have
\be{39a}
\|p^k-\hat p\|_1\,\to\, 0.
\ee

\if{Using the inequality
$$
\|\hat pf_u(\hat w)+\l^kG_u(\hat w)\|_1\le
\| p^kf_u(\hat w)+\l^kG_u(\hat w)\|_1 +\|p_k-\hat p\|_1\| f_u(\hat w)\|_\infty
$$
and conditions (\ref{32a}) and  (\ref{39a}), we get
$\|\hat pf_u(\hat w)+\l^kG_u(\hat w)\|_1\to0.$ This and the condition}\fi

Now, our aim is to find a more detailed representation of the
measure $d\hat\mu.$

\subsection{Representation of the absolutely continuous part of $d\hat\mu$}
\hspace*{1cm} \\
1. Since the sequence $\|\l^k\|_1$ is bounded, then, according to Lemma
\ref{bit}, there exists a sequence of measurable sets $A^k \subset [t_0,t_1]$
such that $\mes A^k\to (t_1-t_0),$ and the sequence $\l_A^k:=\l^k\chi_{A^k}$
is uniformly integrable, hence it contains a weakly convergent (with respect
to $L^\infty)$ subsequence. Without loss of generality we assume that the
sequences  $\l_A^k$ itself weakly converges to some function $\hat \l\in L^1$:
\be{39b}
\l_A^k\, \wky\, \hat \l.
\ee
Since $\l_A^k \ge0$ for all $k,$ we have $\hat \l(t)\ge 0$  {a.e. in} $ [t_0,t_1]$ and
\be{41}
\|\l^k_A\|\, \to\, \|\hat \l\|.
\ee

Further, conditions $\l^k\chi_{M_{\d^k}}=\l^k \ge 0$ (see (\ref{30})),
and $\l_A^k=\l^k\chi_{A^k}$ imply in view of (\ref{Mdk})\, that
$|\l_A^k(t)\,G(\hat w(t))|\le \l_A^k(t)\,\d^k\,$ a.e. on $[t_0,t_1],$
hence $\|\l_A^k\,G(\hat w)\| \to 0.$ The more so,
$\; \l_A^k\,G(\hat w)\to 0\;$  weakly in $L^1.$

On the other hand, (\ref{39b}) implies
$\l_A^k\,G(\hat w)\wky \hat \l\, G(\hat w),$ whence
\be{40}
\hat \l(t)\,G(\hat w(t))\,=\, 0 \q\; \mbox{a.e. in}\;\; [t_0,t_1],
\ee
i.e., the complementary slackness condition in (\ref{16}) holds true.
\ssk

2. Consider more thoroughly condition (\ref{32a}).
Define a function $p_A^k := p^k\chi_{A^k} \in L^\infty.$
Since the set $B^k :=\, [t_0,t_1]\setminus A^k$ has $\mes B^k \to 0,$
and the sequence $\|p^k\|_\infty$ is bounded, we get
$\|p_A^k- p^k\|_1 = \|\,p^k \chi_{B^k}\|_1 \to 0,$ which in view
of (\ref{39a}) yields
\beq{pk1}
\|p_A^k- \hat p\|_1\, \to\, 0.
\eeq
Condition (\ref{32a}) means that
$p^k f_u(\hat w)\,+\, \l^k G_u(\hat w)\; =\, z^k,$ where $\|\,z^k\|_1\, \to 0.$
Multiplying it by $\chi_{A^k},$ we get
$p_A^k\, f_u(\hat w)\,+\, \l^k_A\, G_u(\hat w)\; =\; z^k_A\, :=\,
z^k \chi_{A^k}\,, \q \|\,z^k_A\|_1\, \to 0.$
This and condition (\ref{pk1}) imply
$\|\,\hat p\,f_u(\hat w)\,+\, \l_A^k\,G_u(\hat w)\|_1\,\to\, 0.$
Finally, since $\l^k_A \wky \hat\l,$ we obtain
$\hat p\,f_u(\hat w)\,+\, \hat \l\, G_u(\hat w)\, =\, 0,$
i.e. condition (\ref{20}) of LMP holds true.
\ms

3. Now, introduce the sequence $\l_B^k\,:=\, \l^k\chi_{B^k}\,\in L^1\,.$
Obviously, $\l_B^k\ge 0$ and $\l_A^k\, + \l_B^k\, = \,\l^k.\;$
Note that both $\l_A^k$ and $\l_B^k$ are supported on $M_{\d^k},$ since they
are restrictions of $\l^k$ supported on $M_{\d^k}.$ Therefore, if we narrow
the set $B^k$ to the set $M_{\d^k} \cap B^k,$ the function $\l_B^k$ would not
change. So, we will assume that $B^k \subset M_{\d^k},$
that is $B^k = M_{\d^k} \setminus A^k.\;$ Setting
$$
\dd\mu_A^k\, =\, \l_A^k\,G_x(\hat w)\dd t, \qq
\dd\mu_B^k\, =\, \l_B^k\, G_x(\hat w)\dd t,
$$
we obtain two sequences of measures $\dd\mu_A^k\,$ and $\dd\mu_B^k$ in $C^*.$
Since by (\ref{dmk}) $\dd\mu^k =\l^k\, G_x(\hat w)\dd t,$ we have
$\dd\mu_A^k+\dd\mu_B^k = \dd\mu^k.$ Since $\l_A^k \wky \hat \l,$
we have
\beq{hatdmu1}
\dd\mu_A^k \wst \dd\hat\mu_A\,:=\, \hat \l\, G_x(\hat w)\dd t.
\eeq
Since $\dd\mu^k \wst \dd\hat\mu$ and $\dd\mu_A^k \wst \dd\hat\mu_A,$ there
exists a measure $\dd\hat\mu_B\in C^*$ such that
$$
\dd\mu_B^k \wst \dd\hat\mu_B, \qq \dd\hat\mu_A\,+ \dd\hat\mu_B\,=\dd\hat\mu.
$$
Now we aim to specify the measure $\dd\hat\mu_B\,,$ and this is the main
part of our study.  \ssk

\subsection{Representation of the singular part of $\dd\hat\mu $}
\hspace*{1cm} \\
$1^\circ.$ {We have} $\l_B^k:=\l^k\chi_{B^k}\,,$ where
{$B^k= M_{\d^k} \setminus A^k,$} $\mes B^k\to 0,$ and
\beq{eq37a}
\dd\mu_B^k \,=\,\l_B^k\,G_x(\hat w)\dd t\, \wst \dd\hat\mu_B\,.
\eeq

Since the sequence of norms $\|\l^k\|_1$ is bounded, the sequence of measures
$ \l^k\chi_{B^k}\dd t$ in $C^*$ is also bounded. Therefore, without loss
of generality we assume that there is a measure $\dd\hat\eta\in C^*$ such that
$\dd\hat\eta \ge 0\,$ and
\beq{dmu2eta}
\l_B^k \,\dd t\, =\, \l^k\,\chi_{B^k}\dd t \,\wst\, \dd\hat\eta.
\eeq

Since $\|\l^k\|_1 =\, \|\l_A^k\|_1 + \|\l_B^k\|_1\,,$ conditions (\ref{34a})
and (\ref{41}) imply
\beq{elb1}
\hat\a_0\, +\, \|\hat\l \|_1 \,+\,\|\l_B^k\|_1\,\to 1.
\eeq
Moreover, since $\l_B^k\ge0,$ relation (\ref{dmu2eta}) yields
$\|\l_B^k\|_1\,\to \|\dd\eta\|,$ whence
\beq{elb2}
\hat\a_0\, +\, \|\hat\l \|_1 \,+\,\|\dd\eta\|\,=\, 1,
\eeq
which is equivalent to the nontriviality  condition (\ref{21}).
\ssk

There are two possible cases:\, $\|\dd\eta\|=0$ and $\|\dd\eta\|>0.$
In the first, trivial case, $\|\l_B^k\|_1\,\to 0,$ the more so $\| \dd\mu_B^k\|\to 0,$
then $\dd\hat\mu_B=0,$ i.e. the singular part of $\dd\hat\mu$ does not appear
in the LMP.\, Setting here $\hat s=0$ and $\calR=\O,$ we obtain the costate
equation (\ref{18}) with properties (\ref{17b}) that are trivially satisfied.
\ssk

$2^\circ.$
Consider now the main case, where
\beq{deta-plus}
\|\dd\eta\| \,=\; \lim_k \|\l_B^k\|_1\,=: \,r_B >0.\,
\eeq
Here we will slightly narrow the sets $B^k$ in order to obtain more
properties of $\l^k_B\,.$ To do this, we need the following
\ms

\ble{cut2}\,
Let be given two sequences of functions $a_n\ge0$ and $b_n\ge0$ in
$L^1([t_0,t_1],\R),$ and a sequence of measurable sets $B_n\subset[t_0,t_1]$
of $\mes B_n\to 0$ such that  \vad
$$
\int_{B_n} a_n(t)\dd t\, \to 1, \qq \int_{B_n} b_n(t)\dd t\, \to 0.
$$
Then there is a sequence of measurable sets $E_n\subset B_n$
such that
$$
a_n(t)>0\;\; \mbox{a.e. on}\;\; E_n\,, \q
\int_{E_n} a_n(t)\dd t \to 1,\q \mbox{and} \q
\esssup_{t\in E_n}\, \frac{b_n(t)}{a_n(t)}\,\to\, 0.
$$
\ele

\Proof\, Narrowing if necessary the sets $B_n\,,$ we assume that
$a_n(t)>0$ a.e. on $B_n\,.$ Take any sequence $\w_n\to 0+$ such that
$\dis \int_{B_n} b_n(t)\dd t\, =\, o(\w_n),$ and define a sequence of
sets $E_n =\, \{t\in B_n:\;\, b_n(t) \le \w_n\, a_n(t)\}.\;$ Then
$$
\int_{B_n\setminus E_n} a_n(t)\dd t\, \lee
\frac 1{\w_n}\, \int_{B_n\setminus E_n} b_n(t)\dd t\, \lee
\frac 1{\w_n}\, \int_{B_n} b_n(t)\dd t\; \to\,0,
$$
which gives the required properties. \ctd
\ms

$3^\circ.$ Consider the $L^1-$functions  \vad
\be{47d}
\s^k:=\, p^kf_u(\hat w)+\l^kG_u(\hat w).
\ee
According to (\ref{32a}), $\|\,\s^k\|_1\to0.$ {Then also
$\dis \int_{B^k} (1+|\,\s^k|)\dd t\, \to0.$} By Lemma~\ref{cut2},
there exists a sequence of measurable sets $E^k\subset B^k$
such that $\l_B^k(t)>\,0$ a.e. on $\,E^k,$  \vad
\be{47}
\int_{E^k} \l_B^k\dd t\,\to\, r_B\,>0, \q\; \mbox{and}\q\;
\w^k :=\, \esssup_{E^k}\, \frac{1+|\s^k|}{\l_B^k}\,\to\, 0.
\ee
The first relation here means that
$\dis \int_{B^k\setminus E^k} \l_B^k\dd t\, \to 0,$
therefore (\ref{eq37a}) and (\ref{dmu2eta}) imply
$$
\l^k_B\,G_x(\hat w)\,\chi_{E^k}\dd t\, \wst\, \dd\hat\mu_B\,,\qq
\l^k_B\,\chi_{E^k}\dd t\, \wst\, \dd\hat\eta.
$$
Set  $ \l_E^k\,:=\,\l_B^k\chi_{E^k}\,=\,\l^k\chi_{E^k}\,.\;$
Then $\l^k(t)=\l_E^k(t)$  a.e. on $\,E^k,$
\be{47m}
 \l^k_E\dd t \,\wst\, \dd\hat\eta, \qq
 \l_E^k\,G_x(\hat w)\dd t \, \wst\, \dd \hat \mu_B\,,
\ee
so the ''narrowed'' sequence $ \l^k_E$ has the same limit properties
as the original $\l^k_B$ does.
\ssk

Since $E^k \subset M_{\d^k},$ relations (\ref{Mdk}) imply
\beq{46f}
-\d^k \le G(\hat w(t))\le 0\q\; \mbox{a.e. on}\q E^k.
\eeq
Moreover, in view of definition (\ref{47d}), for all $k$
$$
G_u(\hat w)\,=\, \frac{\s^k-p^kf_u(\hat w)}{\l^k} \qq
\mbox{a.e. on} \q E^k.
$$
The second condition in (\ref{47}) and the boundedness of the
sequence $\|p^k\|_\infty$ imply
$$
\varepsilon^k :=\, \esssup_{E^k} |\,G_u(\hat w)|\,=\;
\esssup_{E^k}\frac{\left|\s^k-p^kf_u(\hat w)\right|}{\l^k}\,\to\, 0,
$$
whence  \vad
\beq{47f}
|\,G_u(\hat w(t))|\lee \varepsilon^k \q\; \mbox{a.e. on}\q E^k.
\eeq \ssk

$4^\circ.$ We will need the following constructions. \,
Recall that {we introduced the set of phase points}
${\calN(G)}:=\{(x,u)\in\R^{n+m}:\; G(x,u)=0, \;\, G_u(x,u)=0 \}\,$
{and assumed that this set is nonempty.}

Now, for any $\d>0$ and  $\e>0,$ introduce {its extension} up to $\delta,\e$:
$$
\calN_{\d,\e}(G)\,:=\; \{ (x,u)\in\R^{n+m}\,: \;\; -\d\le G(x,u)\le 0,
\q  |\,G_u(x,u)|\le\e \,\}.
$$
Obviously it is closed, and
$\calN(G)\, =\, \bigcap_{\d>0,\;\, \e>0} \calN_{\d,\e}(G).\;$
\ssk

By analogy with the mapping $S(x,u),$ for any $\d>0$ and $\e>0,$
define a set-valued mapping  \vad
$$
(x,u)\in \R^{n+m}\;\, \To\;\, S_{\d,\e}(x,u) \,\subset \R^{n*}:
$$
\begin{itemize}
   \item[(i)] if    $(x,u)\in \calN_{\d,\e}(G)\,,\;$
   then $S_{\d,\e}(x,u)\,=\,\{G_x(x,u)\},$ \ssk
   \item[(ii)] if $(x,u)\notin \calN_{\d,\e}(G)\,,\;$
   then $S_{\d,\e}(x,u)\, =\, \O.$
\end{itemize} \ssk

\noi
Obviously, this mapping is compact-valued, upper semicontinuous, and
\beq{Sdexu}
\bigcap_{\d>0,\;\, \e>0} S_{\d,\e}(x,u)\, =\, S(x,u) \qq
\mbox{for all}\;\, (x,u).
\eeq
For any nonempty set $M\subset \R^{n+m},$ define
$$
S_{\d,\e}(M):= \bigcup_{(x,u)\in M} S_{\d,\e}(x,u).
$$
and $S_{\d,\e}(\O)=\O.$ Clearly, if $M$ is compact, the set
$S_{\d,\e}(M)$ is compact as well.\,
Note that for any  $M\subset \R^{n+m}$ we have
$
S_{\d,\e}(M)\; =\; S_{\d,\e}(M\cap \calN_{\d,\e}(G)).
$

Now, consider the reference process $\hw(t) = (\hx(t),\hu(t)).\,$
We will assume that the corresponding set $\calD,$ defined in (\ref{setD}),
is nonempty,  i.e. there exists a point $ t_*\in[t_0,t_1]$ such that
$\; \clm (\hw)(t_*)\, \cap\, \calN(G)\ne\O.$
\ssk

Since the set $\clm(\hat u)$ is compact, the set $S_{\d,\e}(\hat x(t),\,\clm(\hat u)(t))$
is also compact for any $t$ and upper semicontinuous in~$t.$ \ssk

For any points $\t_0 <\t_1\,$ in $[t_0,t_1],$ define a set
\beq{Q-de}
Q_{\d,\e}[\t_0,\t_1]\, :=\;
\bigcup_{t\in[\t_0,\t_1]} S_{\d,\e}(\hat x(t),\,\clm(\hat u)(t)).
\eeq
By the above argument, the right hand side here is a compact set.
Moreover, relation (\ref{Sdexu}) implies that  \vad
\beq{Qxu-cap}
Q[\t_0,\t_1]\,:=\; \bigcap_{\d>0,\;\, \e>0}\, Q_{\d,\e}[\t_0,\t_1]\, =\;
\bigcup_{\t\in[\t_0,\t_1]} S(\hat x(\t),\,\clm(\hat u)(\t)),
\eeq
and this set is also compact.\,

Finally, for any $t_*\,,$ if $\t_0 \to t_*-0$ and $\t_1\to t_*+0,$
then obviously
\beq{Qtau01}
Q[\t_0,\t_1]\; \to\; S(\hat x(t_*),\,\clm(\hat u)(t_*)).
\eeq
in the Hausdorf sense.\,
\ms

$5^\circ.$ Now we can describe the relationship between the measures
$\dd \hat \mu_B$ and $\dd\hat\eta.$
{Recall that $\|\dd\hat\eta\| >0$ by (\ref{deta-plus}). }
\ms

\ble{deta} \q 
The measure $\dd\hat\mu_B$ admits a representation
\be{47m2}
\dd\hat\mu_B\; =\; \hat s(t)\dd\hat\eta \q\; 
\ee
with some $\dd\hat\eta$-measurable essentially bounded function
$\hat s: [t_0,t_1]\to \R^{n*},$ and there is a set $\calR\subset \calD$ of full
$\dd\hat\eta$-measure (i.e., $\int_\calR\dd\hat \eta = \int_{[t_0,t_1]}\dd\hat \eta$)\,
such that
\beq{scoclm}
\hat s(t)\in\conv\, S\big(\hat x(t),\,\clm(\hat u)(t)\big) \q
\mbox{for all}\q t\in\calR.
\eeq
\ele

\Proof\, a) In view of (\ref{47m}), 
$|\dd\hat\mu_B| \le M\dd\hat\eta,$ where $M = \|G_x(\hat w)\|_\infty\,.$
Hence, the measure $\dd\hat\mu_B$ is absolutely continuous with
respect to $\dd\hat\eta$, and therefore, by the Radon-Nikodym  theorem
it admits representation (\ref{47m2}), where $\hat s(t)$ is
a $\dd\hat\eta$-measurable function taking values in $\R^{n*}$ and
satisfying $|\hat s(t)|\le M$ a.e. in $\dd\hat\eta.$ \\[4pt]
\hsp b) Let us prove inclusion (\ref{scoclm}) with some $\calR\subset \calD$
of full	$\dd\hat\eta$-measure.\,
Fix a point $t_*\in[t_0,t_1]$ with the following properties: \ssk

\begin{itemize}
    \item[(i)] if $t_*\in[\tau_0,\tau_1]\subset [t_0,t_1],$  $\tau_0<\tau_1,$
then $\int_{[\tau_0,\tau_1] } \dd\hat\eta>0,$

    \item[(ii)] if $t_*\in[\tau_0,\tau_1] \subset [t_0,t_1],$ $\tau_0<\tau_1\,,$
    $\tau_0\to t_*,$ $\tau_1\to t_*,$ then\footnote{\,If $\t_0 =t_*\,$
    we do not need to tend $\tau_0\to t_*,$ so only tend $\tau_1\to t_*\,.$
    The same concerns $\t_1$.}
\be{48a}
\hat s_{[\tau_0,\tau_1] }\,:=\; \frac{\int_{[\tau_0,\tau_1]} \dd\hat\mu_B}
{\int_{[\tau_0,\tau_1] } \dd\hat\eta}\;\, \to\, \hat s(t_*),
\ee
\end{itemize}
As is known, the set of such points $t_*$ has a full $\dd\hat\eta$-measure
in $[t_0,t_1]$ (since it includes the Lebesgue points of the function
$\hat s$ with respect to the measure $\dd\hat\eta).$ Denote this set by
$\calR.\,$ We have $\int_\calR\dd\hat\eta\,=\,\int_{[t_0,t_1]}\dd\hat\eta.\,$
\ssk

Take any $[\t_0,\,\t_1]$ containing $t_*\,.$ Then
$\int_{[\tau_0,\tau_1] }\,  \l^k_E\,\dd t>0$ for all sufficiently
large~$k.$ (Otherwise $\l^k_E=0$ a.e. in $[\t_0,\,\t_1]$ for a subsequence
$k\to \infty,$ which implies that also $\dd \hat\eta =0$ in $[\t_0,\,\t_1],$
a contradiction with the choice of $t_*).$
\ssk

Therefore, we can define a row-vector
\be{49s}
s^k_{[\tau_0,\tau_1]}\,:= \;\frac{\int_{[\t_0,\t_1]}\,\l^k_E\,G_x(\hw)\,\dd t}
{\int_{[\tau_0,\tau_1] }\,  \l^k_E\,\dd t }\;.
\ee

Let $\Theta$ be the set of continuity of the measures $\dd\hat\mu_B$ and
$\dd\hat\eta,$ i.e., the set of all those $t,$ which are not atoms neither
of $\dd\hat\mu_B$ nor of $\dd\hat\eta.$ Note that $\Theta$ is dense in
$[t_0,t_1].\,$ If $\t_0,\, \t_1 \in \Theta,$ then
\be{49}
s^k_{[\tau_0,\tau_1] }\,\to\, \hat s_{[\tau_0,\tau_1] }\q\;
\mbox{as}\q k\to\infty,
\ee
since both the numerator and denominator tend to the corresponding limits.
\ssk

Notice that the right hand side of (\ref{49s}) is a ``convex combination''
of the vectors $G_x(\hw(t)),$ in its continuous version.
\ssk

c) In view of (\ref{46f}) and (\ref{47f}), we have
$\,\hw(t) \in \calN_{\d^k\e^k}(G)\;$ a.e. on $E^k,\,$ and so
$$
G_x(\hw(t))\,\in\,  S_{\d^k \e^k}(\hat x(t),\hat u(t)) \q\mbox{a.e. on} \q E^k.
$$
Since $\;\hat u(t)\in \clm(\hat u)(t)\;$ a.e. on $[t_0,t_1],\,$
we get
\beq{GxS}
G_x(\hat w(t))\,\in\, S_{\d^k\e^k}(\hat x(t),\,\clm(\hat u)(t))
\q\mbox{a.e. on} \q E^k.
\eeq

Recall that $t_* \in \calR.\,$ Take any $[\t_0,\,\t_1]$ containing $t_*\,.$
The last inclusion and definition (\ref{Q-de}) imply that
for almost all $t\in E^k \cap [\t_0,\,\t_1]$
$$
G_x(\hat w(t))\,\in\,
\bigcup_{\t\in[\t_0,\t_1]} S_{\d^k\e^k}(\hat x(\t),\,\clm(\hat u)(\t))
\;=\; Q_{\d^k\e^k}[\t_0,\t_1]\,,
$$
and the more so,\, for almost all $t\in E^k \cap [\t_0,\,\t_1]$
$$
G_x(\hat w(t))\,\in\, \conv\, Q_{\d^k\e^k}[\t_0,\t_1].
$$
By the Caratheodory theorem, the right hand side here is a convex compact
set. Then, the definition (\ref{49s}) gives (since $\l^k_E$ is supported
on $E^k):$
$$
s^k_{[\tau_0,\tau_1] }\,\in\, \conv\, Q_{\d^k\e^k}[\t_0,\t_1]\,.
$$
Now, assume that $\t_0,\,\t_1\in \Theta,\,$ $\t_0 <t_*< \t_1\,.$
Taking the limit as $k\to\infty\,$ in view of (\ref{49})
and (\ref{Qxu-cap}), we get  \vad
$$
\hat s_{[\tau_0,\tau_1] }\,\in \, \conv\, Q(x,u)[\t_0,\t_1]\,.
$$
Finally, taking the limit as $\t_0\to t_*\,,\,$  $\t_1\to t_*\,$ along $\Theta\,$
in view of (\ref{48a}) and (\ref{Qtau01}), we obtain  \vad
\be{49ss}
\hat s(t_*)\, \in\, \conv S(\hat x(t_*),\,\clm(\hat u)(t_*)).
\ee
Thus, the set $ S(\hat x(t_*),\,\clm(\hat u)(t_*))$ is nonempty, which
by the definition (\ref{setD}) means that $t_*\in \calD.$ Since the point
$t_*\in\calR$ is arbitrary, it follows that $\calR \subset \calD.$
Consequently, inclusion (\ref{scoclm}) holds.
\ssk

d) For $t\in [t_0,t_1]\setminus\calR$ we can redefine (if necessary) $\hat s(t)$
by zero, without violating the conditions of LMP.\,  The lemma is proved.
\ctd \ms

Thus, in view of this lemma and (\ref{hatdmu1}), the adjoint equation
(\ref{38}) takes the form
\be{50}
-\dd \hat p\, =\; \hat pf_x(\hat w)+\dd\hat \mu_A+ \dd\hat \mu_B\; =\;
\hat pf_x(\hat w)+\hat \l\, G_x(\hat w)\dd t+\hat s\dd\hat\eta,
\ee
i.e. condition (\ref{18}) holds true.  \ms

Thus, all conditions (\ref{15})--(\ref{21}) of the local minimum principle
are satisfied.\\ Theorem~\ref{thmp}\, is completely proved.
\bs


\end{document}